\pdfoutput=1
\RequirePackage{ifpdf}
\ifpdf
\documentclass[pdftex]{sigma}
\else
\documentclass{sigma}
\fi

\usepackage{amsmath, mathtools}
\usepackage{lmodern}
\usepackage{microtype}
\usepackage{booktabs}
\usepackage{enumerate}
\usepackage{amscd}
\usepackage{pb-diagram}
\usepackage[all, cmtip]{xy}
\usepackage{esint}
\usepackage{hypcap}
\usepackage{tikz}
\usetikzlibrary{calc,decorations.pathmorphing,shapes}

\numberwithin{equation}{section}

\newcommand{\C}{\mathbb{C}}

\newcommand{\divv}{\operatorname{div}} % \div is a symbol; avoid redefining it

\newcommand{\Spec}{\operatorname{Spec}}
\newcommand{\Hol}{\operatorname{Hol}}
\newcommand{\Hom}{\operatorname{Hom}}
\newcommand{\Cl}{\operatorname{Cl}}
\newcommand{\Supp}{\operatorname{Supp}}
\newcommand{\Pic}{\operatorname{Pic}}

\newtheorem{Theorem}{Theorem}[section]
\newtheorem{TibaTheorem}[Theorem]{Theorem}
\newtheorem{Proposition}[Theorem]{Proposition}
\newtheorem{Corollary}[Theorem]{Corollary}
\newtheorem{Lemma}[Theorem]{Lemma}
\theoremstyle{definition}
\newtheorem{Definition}[Theorem]{Definition}
\newtheorem{Remark}[Theorem]{Remark}

\newcounter{sarrow}

\begin{document}

\allowdisplaybreaks

\newcommand{\arXivNumber}{2111.02192}

\renewcommand{\PaperNumber}{119}

\FirstPageHeading

\ShortArticleName{Hyperbolic Embeddings in Toric Geometry: Effectivity and Deformation Stability}

\AuthorNameForHeading{G.~Cho and J.~Yoo}

\ArticleName{Hyperbolic Embeddings in Toric Geometry:\\ Effectivity and Deformation Stability}

\Author{Gunhee CHO~$^{\rm a}$ and Jaewon YOO~$^{\rm b}$}
\Address{$^{\rm a)}$~Department of Mathematics, Texas State University,\\
\hphantom{$^{\rm a)}$}~601 University Dr, San Marcos, TX~78666, USA}
\EmailD{\mail{wvx17@txstate.edu}}

\Address{$^{\rm b)}$~Department of Mathematics, Pohang University of Science and Technology,\\
\hphantom{$^{\rm b)}$}~Jigok-ro 20, Pohang-Si, Gyeongsangbuk-do 37665, South Korea}
\EmailD{\mail{yooj1215@postech.ac.kr}}

\ArticleDates{Received February 24, 2026, in final form November 18, 2026; Published online November 22, 2026}

\Abstract{We study the deformation behavior of Kobayashi hyperbolic embeddings for complements of divisors in projective toric varieties. In the toric setting, the Zariski closure of an entire curve in a~torus complement is a~translate of an algebraic subtorus contained in that complement. This translates hyperbolicity questions into combinatorial conditions on lattice-point configurations of Newton polytopes. Building on a~theorem of Tiba, which guarantees hyperbolic embedding for a~	general divisor under suitable facewise lattice conditions, we develop an effective refinement of his argument. We construct an explicit constructible bad locus in the coefficient parameter space, characterized by the presence of translated subtori in the divisor 	or its complement on a~torus stratum, and use its proper Zariski closure as an exceptional locus. This description makes the exceptional set amenable to computation. Using this effectivity, we prove a~deformation stability result: along any algebraic one-parameter family of divisors whose initial member avoids the exceptional locus, hyperbolic embedding fails for at most finitely many parameters. Under strengthened lattice-point hypotheses, we further exhibit distinguished one-parameter families for which hyperbolic embedding persists without any exceptional parameters. We also analyze deformations arising from diagonal torus reparametrizations. If the initial divisor lies outside $\mathcal Y$, the inverse image of $\mathcal Y$ under the similarity family is empty on the nondegenerate parameter domain, while the smallest possible degeneracy set of a~rational similarity presentation is controlled by a~multiplicative rank invariant of the scaling vector. Finally, we illustrate the theory through examples in projective spaces and Hirzebruch surfaces. For cubic surfaces in $\mathbb P^3$, a~sparse-resultant calculation shows that the reduced divisorial part of the global exceptional locus has degree~$240$.}
	
\Keywords{Kobayashi metric; Kobayashi hyperbolicity; hyperbolic embeddings in toric varieties}

\Classification{32Q45; 14M25}
	
\section{Introduction and results}\label{sec:intro-results}

\subsection{Deformation and (non)hyperbolicity in the noncompact world}

Deformation theory has long played a~decisive role in the study of
(non)hyperbolicity in complex geometry. A~precise propagation result for the Kobayashi pseudodistance was proved by
Kamenova--Lu--Verbitsky: if this pseudodistance vanishes identically for one
complex structure, then it vanishes for every ergodic complex structure in
the same deformation class~\cite[Theorem~2.1]{MR3263959}.
They deduced, in particular, that all complex K3 surfaces have identically
vanishing Kobayashi pseudodistance~\cite[Corollary~2.2]{MR3263959}.
More recently, one of the authors together with D.R.~Morrison proved
non-measure-hyperbolicity for every complex K3 surface $S$ and every Hilbert
scheme of points $\operatorname{Hilb}^n(S)$, $n\ge1$
\cite[Theorem~2 and Corollary~3]{ChoMorrison2025}.
These results illustrate the role of deformation methods in proving
vanishing statements.

The situation is markedly different for hyperbolicity.
For compact complex manifolds, Kobayashi hyperbolicity is open under deformation.
For general deformations, upper semicontinuity of the infinitesimal
Kobayashi--Royden pseudometric is available~\cite{WrightMarcus77};
see also~\cite{KalkaM76} for the compact hyperbolic case.
For divisor complements, however, hyperbolicity need not persist under
small deformations. Zaidenberg gives an explicit pencil deforming five
lines in general position in $\mathbb P^2$ to smooth quintics whose
complements contain an affine line
\cite[Proposition~3.2 and Example~3.2]{Zaidenberg1989}.
Thus, deformation theory in the noncompact world is a~genuinely delicate
operation, and stability of hyperbolicity cannot be taken for granted.

This leads naturally to the central question addressed in this paper:
\emph{Under what conditions does Kobayashi hyperbolicity of a~noncompact
		complement persist under algebraic deformation of divisors?}

\subsection{Toric geometry and combinatorial control of entire curves}

We study this question in the toric setting, where analytic properties of
divisor complements admit a~precise translation into combinatorial data.
Throughout, let $X=X_\Sigma$ be a~projective toric variety of dimension
$r\ge1$ with dense torus $T_N\cong\bigl(\mathbb{C}^*\bigr)^r$ and character lattice
$M\cong\mathbb{Z}^r$, and let~$D$ be a~torus-invariant ample Cartier divisor.
Its polytope $P_D\subset M_{\mathbb R}$ is a~full-dimensional lattice
polytope with normal fan $\Sigma$.
We consider divisors in $|D|$ represented by Laurent polynomials supported
in $P_D\cap M$; the actual Newton polytope of a~member may be smaller when
some coefficients vanish.

In this framework, the Zariski closure of an entire curve in a~torus
complement is a~translate of an algebraic subtorus contained in that
complement; see Theorem~\ref{thm:DethloffLu}.
As a~consequence, the exclusion of entire curves can be detected via
lattice-point configurations along faces of the Newton polytope. Combined
with Tiba's boundary criterion on all torus strata, this yields Kobayashi
hyperbolic embeddings.
This philosophy was made precise by Tiba, who proved that if a~finite
support $S\subset P_D\cap M$ satisfies suitable facewise combinatorial
conditions, then the complement of a~\emph{general} divisor in the
corresponding linear system is Kobayashi hyperbolically embedded into $X$
(see Theorem~\ref{thm:Tiba13MainThm}).

Tiba's proof already uses layer polynomials, incidence dimension estimates,
and a~finite collection of hyperplane directions to obtain a~proper
exceptional subset. We use this construction as our starting point.
Our emphasis is on distinguishing its constructible incidence images from
their Zariski closures, making their coefficient conditions usable in
specific families, and carrying out explicit elimination calculations.

\subsection{Main contributions and structure of the paper}

We give a~coefficient-wise reformulation of Tiba's argument and apply it
to special deformations and explicit examples. The generic hyperbolic
embedding criterion remains Tiba's theorem; the finite-exception statement
below is an algebraic consequence of a~proper closed exceptional locus.

Our first contribution is an effective reformulation of Tiba's argument.
Under the same facewise lattice hypotheses, we construct a~constructible bad
locus and a~proper Zariski closed exceptional locus
$\mathcal Y \subset |D|$
in the coefficient parameter space. The bad locus is defined by concrete
facewise degeneration conditions that detect the presence of translated
subtori, and $\mathcal Y$ is its reduced Zariski closure (together with the
boundary of the relevant coefficient chart).
This description makes the exceptional set amenable to explicit
computation, a~feature illustrated by several concrete examples.

The construction has the following deformation consequence. Along an
algebraic one-parameter family whose image is not contained in
$\mathcal Y$, hyperbolic embedding holds outside a~finite set of parameters.
The hypothesis on the initial member in Theorem~\ref{thm:1} ensures this
noncontainment.

\begin{Theorem}[deformation stability with finite exceptions]
	\label{thm:1}
	Assume that, for every positive-dimensional face $\tau$ of $P_D$ and
	every $H\in\mathcal H_{\tau\cap M}$,
	\[
	\bigl|\phi_H((\tau\cap M)-x)\bigr|\ge \dim\tau+1
	\]
	for some $x\in\tau\cap M$. $($Equivalently, the complete support
	$P_D\cap M$ satisfies the facewise hypotheses of
	Theorem~{\rm\ref{thm:Tiba13MainThm}}.$)$
	Let $D_0\in\mathfrak d_D\setminus\mathcal Y$, where $\mathcal Y\subset |D|$
	is the proper Zariski closed exceptional locus constructed in
	Section~{\rm\ref{sec:proof-main}}.
	Let $\gamma\colon \mathbb P^1\to |D|$ be a~morphism with $\gamma(0)=D_0$.
	Then for all but finitely many $t\in\mathbb P^1$, the complement
	\[
	T_N\setminus\Supp(D_t\cap T_N),\qquad D_t:=\gamma(t),
	\]
	is Kobayashi hyperbolically embedded in $X$.
\end{Theorem}

The same conclusion holds for a~morphism $\gamma\colon C\to|D|$ from an
irreducible complex algebraic curve $C$, provided some member lies outside
$\mathcal Y$. Here and below algebraic curves are of finite type over
$\mathbb C$. If $C$ is reducible, the hypothesis is required on each
irreducible component to which the conclusion is applied; see
Section~\ref{subsec:proof-finite-exceptions} for the proof.

Finally, under strengthened lattice-point conditions, we show that
exceptional parameters can be eliminated altogether for certain
distinguished affine one-parameter families.

\begin{Theorem}[exception-free deformations]\label{thm:2}
	Let $S\subset P_D\cap M$ be finite and meet every facet of $P_D$
	$($equivalently, $\mathfrak d_S\ne\varnothing${\rm)}.
	Assume that for every
	positive-dimensional face $\tau$ of $P_D$, the set~${S_\tau:=\tau\cap S}$ satisfies $\operatorname{aff}(S_\tau)=\operatorname{aff}(\tau)$ and $
	|\phi_H(S_\tau-x)|\ge \dim\tau+2$ for every $H\in\mathcal H_{S_\tau}$ and~${x\in S_\tau}$.
	Choose a~representative $a=(a_I)_{I\in S}$ whose projective class
	$[a]\in\mathbb P^{|S|-1}$ lies outside the enlarged facewise
	exceptional locus $\mathcal Y'_S$ of
	Definition~{\rm\ref{def:global-enlarged-locus}}, and fix $I_0\in S$. Then
	\[
	D_t=\divv\biggl(tz^{I_0}+\sum_{I\in S\setminus\{I_0\}}a_Iz^I\biggr),
	\qquad t\in\mathbb C,
	\]
	where $D_t$ denotes the divisor of the corresponding global section of
	$\mathcal O_X(D)$. For every $t\in\mathbb C$, the inclusion
	$T_N\setminus\Supp(D_t\cap T_N)\hookrightarrow X$ is a~Kobayashi
	hyperbolic embedding.
\end{Theorem}

The original divisor occurs at $t=a_{I_0}$, whereas $t=0$ deletes that
coefficient. The assertion does not include the member at infinity of the
projective closure of this family.
 A~key feature of our approach is that the bad loci are described by explicit
incidence conditions. We give examples in $\mathbb P^2$, $\mathbb P^3$, and
Hirzebruch surfaces; in the cubic-surface case we compute the reduced
divisorial part of the global exceptional locus by sparse resultants and obtain
degree $240$.

\textbf{Structure of the paper.}
Section~\ref{sec:toric-geometry} recalls basic toric notation and fixes the
Newton polytopes and combinatorial setup used throughout.
In Section~\ref{sec:proof-main}, we give an effective proof of Tiba's theorem
by constructing the exceptional locus $\mathcal Y$ in coefficient space and
relating its defining conditions to translated subtori via the theorem of
Dethloff--Lu recalled below.
Section~\ref{sec:deformation-stability} is devoted to deformation stability and
contains the proofs of Theorems~\ref{thm:1} and~\ref{thm:2}, including the refined
enlarged exceptional locus tailored to one-monomial borderline degenerations.
Section~\ref{sec:similarity-exceptions} studies a~distinguished class of
one-parameter deformations arising from diagonal torus reparametrizations
(\emph{similarity}). For an initial member outside $\mathcal Y$, the inverse image of
$\mathcal Y$ under the similarity family is empty on the nondegenerate domain; the section relates the minimal number of
degenerate parameters in a~rational presentation to the multiplicative rank
of the scaling vector $\lambda$.
Finally, Section~\ref{sec:examples} illustrates the theory by explicit computations
of exceptional loci and deformation behavior in concrete toric models, including
projective spaces and Hirzebruch surfaces; it also contains the
sparse-resultant computation of the degree-$240$ divisorial exceptional locus
for cubic surfaces in~$\mathbb P^3$.

\section{Kobayashi hyperbolicity and hyperbolic embeddings}

Throughout this paper, we work over the field $\mathbb{C}$.
We briefly recall the notions of the Kobayashi--Royden metric,
Kobayashi hyperbolicity, and hyperbolic embedding, emphasizing
their role in deformation problems on \emph{noncompact} complex manifolds.

\subsection{Kobayashi--Royden metric and Kobayashi distance}

Let $X$ be a~complex manifold.
The \emph{Kobayashi--Royden $($pseudo-$)$metric} $KR_X(p;v)$ at a~point $p\in X$ in the direction $v\in T_pX$ is defined by
\[
KR_X(p;v)
=
\inf\left\{
\frac{1}{\alpha} \mid
\alpha>0,\,
f\in\Hol(\mathbb{D},X),\,
f(0)=p,\,
f'(0)=\alpha v
\right\},
\]
where $\mathbb{D}\subset\mathbb{C}$ denotes the unit disk.
The pseudometric measures derivatives of holomorphic disks centered at
$p$ and tangent to $v$; these maps need not be injective
(see, for example, \cite{GunheeChoJunqingQian20}).
In general, $KR_X$ is only a~pseudometric: it may vanish for
nonzero tangent vectors.

The associated \emph{Kobayashi pseudodistance} $d_X$ on $X$ is defined by
\[
d_X(x,y)
=
\inf
\left\{
\sum_{i=1}^n \rho_{\mathbb{D}}(a_i,b_i)
\right\},
\]
where $x=p_0,\dots,p_n=y$, and for each $i$ there exists a~holomorphic map
$f_i\in\Hol(\mathbb{D},X)$ such that~${f_i(a_i)=p_{i-1}}$ and $f_i(b_i)=p_i$.
Here $\rho_{\mathbb{D}}$ denotes the Poincar\'e distance on $\mathbb{D}$.
It is known that $d_X$ coincides with the length distance induced by
the Kobayashi--Royden metric~\cite[Theorem~3.1.15]{KobayashiShoshichi98}.

\subsection{Kobayashi hyperbolicity}
 A~complex manifold $X$ is said to be \emph{Kobayashi hyperbolic}
if the pseudodistance $d_X$ is a~genuine distance, i.e., separates points.
Kobayashi hyperbolicity implies that $X$ admits no nonconstant entire
holomorphic curve $f\colon \mathbb{C}\to X$. The converse is the Brody criterion
when $X$ is compact, but it is false for arbitrary noncompact complex
manifolds. We therefore use the term \emph{Brody hyperbolic} for the absence
of nonconstant entire curves and keep it distinct from Kobayashi
hyperbolicity in the noncompact setting.

For compact complex manifolds, Kobayashi hyperbolicity is open under deformation.
For general deformations, the infinitesimal Kobayashi--Royden
pseudometric is upper semicontinuous~\cite{WrightMarcus77};
see also~\cite{KalkaM76} in the compact hyperbolic case.
This statement alone does not ensure deformation openness of
hyperbolicity for noncompact fibers.

Zaidenberg's pencil of quintics gives an explicit example in which a~small deformation of five lines in general position destroys
hyperbolicity of the complement
\cite[Proposition~3.2 and Example~3.2]{Zaidenberg1989}.
Consequently, in noncompact geometry one must distinguish
between intrinsic hyperbolicity and the behavior of the manifold
inside a~compact ambient space.

\subsection{Hyperbolic embeddings}
\label{subsec:hyperbolic-embeddings}

Let $X$ be a~complex manifold and $Y\subset X$ a~complex submanifold,
with closure $\overline{Y}\subset X$.
Following~\cite[Section~3.3]{KobayashiShoshichi98},
we say that $Y$ is \emph{Kobayashi hyperbolically embedded} in $X$
if for any two distinct points $p,q\in\overline{Y}$ there exist
neighborhoods $U_p$ of $p$ and $U_q$ of $q$ such that
\[
d_Y(U_p\cap Y,U_q\cap Y)>0.
\]
The same neighborhood criterion defines hyperbolic embedding when the
ambient space $X$ is a~possibly singular complex space, as in the toric
applications below.

If $Y\hookrightarrow X$ is a~hyperbolic embedding, then $Y$ itself
is Kobayashi hyperbolic.
The converse, however, is false in general
\cite[Example~3.3.10]{KobayashiShoshichi98}.
From the perspective of deformation theory, hyperbolic embedding
is the correct notion: it controls not only entire curves inside $Y$,
but also their possible degenerations toward the boundary
$\overline{Y}\setminus Y$; see also~\cite{Zaidenberg1986}.

Our guiding principle is the following. In this paper, we study deformation
problems for complements
$T_N\setminus D$ of divisors in toric varieties.
Our central philosophy is that \emph{hyperbolic embedding} provides a~robust replacement for Kobayashi hyperbolicity in noncompact settings,
allowing effective control under deformation.

% ======================================================
\section{Toric geometry and hyperbolic embedding}
\label{sec:toric-geometry}

In this section, we fix notation and record the toric-geometric framework
used throughout the paper.
Our purpose is to make precise how combinatorial data of lattice polytopes
control the behavior of entire curves in complements of divisors.
In particular, we recall the hyperplane-projection formalism underlying
Tiba's generic hyperbolicity criterion, which will serve as a~baseline for
the deformation-theoretic refinements developed in later sections.

Let $X=X_\Sigma$ be a~projective toric variety of dimension $r\ge1$, and let
$T_N \cong \bigl(\mathbb{C}^*\bigr)^r$ denote its dense algebraic torus.
Let
$
M \cong \mathbb{Z}^r$, $ N \cong \mathbb{Z}^r$
be the character lattice and the lattice of one-parameter subgroups,
respectively, equipped with the canonical pairing
$\langle\cdot,\cdot\rangle \colon M \times N \to \mathbb{Z}$.
Set
$
M_{\mathbb{R}} := M \otimes_{\mathbb{Z}} \mathbb{R}$, $
N_{\mathbb{R}} := N \otimes_{\mathbb{Z}} \mathbb{R}$.

 \textbf{Polytopes and divisors.}
As fixed in the introduction, $D$ is a~torus-invariant ample Cartier
divisor on $X$. Write
\smash{$
D = \sum_{\rho \in \Sigma(1)} a_\rho D_\rho$},
where $\Sigma(1)$ denotes the set of rays of $\Sigma$ and $D_\rho$ is the
corresponding invariant prime divisor.
Let~${m_\rho \in N}$ be the primitive generator of $\rho$.
The ampleness and Cartier assumptions imply that $D$ determines a~full-dimensional lattice polytope
\[
P_D
=
\bigcap_{\rho \in \Sigma(1)}
\bigl\{
u \in M_{\mathbb{R}}
\mid
\langle u, m_\rho \rangle \ge -a_\rho
\bigr\}
\subset M_{\mathbb{R}},
\]
whose normal fan coincides with $\Sigma$; see~\cite[Chapter~2]{FultonToric}.

Via the Cox construction~\cite[Section~5.1]{CoxToricVarieties},
we may realize
\[
X
\cong
\bigl(
\Spec \mathbb{C}[x_\rho \mid \rho \in \Sigma(1)] \setminus Z(\Sigma)
\bigr)
/\!\!/ G(\Sigma),
\]
where $Z(\Sigma)$ is the exceptional set and
$G(\Sigma)=\Hom_{\mathbb{Z}}\bigl(\Cl(X),\mathbb{C}^*\bigr)$.

Each lattice point $I \in P_D \cap M$ corresponds to a~monomial
\[
x^{\langle I, P_D \rangle}
=
\prod_{\rho \in \Sigma(1)}
x_\rho^{\langle I, m_\rho \rangle + a_\rho},
\]
whose divisor is linearly equivalent to $D$.

 \textbf{Finite supports and linear systems.}
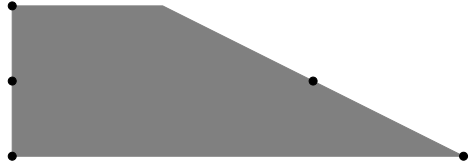
Let $S \subset P_D \cap M$ be a~finite set of lattice points.
For a~nonzero vector $a=(a_I)_{I\in S}$, write
\[
s_a:=\sum_{I\in S}a_Ix^{\langle I,P_D\rangle},\qquad
D_X(a):=\operatorname{div}_X(s_a),\qquad
D_T(a):=D_X(a)|_{T_N}.
\]
We denote the effective Weil-divisor closure of $D_T(a)$ in $X$ by
\smash{$D_{\mathrm{cl}}(a):=\overline{D_T(a)}^{ X}$}, retaining the
multiplicities of its components. The difference
$D_X(a)-D_{\mathrm{cl}}(a)$ consists precisely of the torus-invariant
components of $D_X(a)$. We also write $D(a)=D_X(a)$ when referring to a~member of the coefficient linear system.
We define the projective linear system
\[
\mathfrak{d}_S^*
=
\biggl\{
\divv \biggl(\sum_{I\in S} a_I x^{\langle I,P_D\rangle}\biggr)
\mid
(a_I)_{I\in S}\in\mathbb{P}^{|S|-1}
\biggr\}.
\]
Inside $\mathfrak d_S^*$, let
\[
\mathfrak d_S
:=
\bigl\{D(a)\in\mathfrak d_S^*
\mid
(a_I)_{I\in\tau\cap S}\ne 0
\text{ for every facet }\tau\subset P_D
\bigr\}.
\]
Thus $\mathfrak d_S$ consists precisely of the members having no
torus-invariant irreducible component. Equivalently,
$D_X(a)=D_{\mathrm{cl}}(a)$ for these members. Restricting a~member of
$\mathfrak d_S$ to $T_N$ yields the Laurent polynomial
\[
\sum_{I\in S} a_I z^I,
\qquad
z^I = z_1^{i_1}\cdots z_r^{i_r},
\]
whose zero divisor determines the member on the dense torus. In particular,
$\mathfrak d_S$ is nonempty if and only if every facet $\tau$ of $P_D$
satisfies $\tau\cap S\ne\varnothing$. We abbreviate
$\mathfrak d_D:=\mathfrak d_{P_D\cap M}$.
For the global facewise constructions we require
$\mathfrak d_S\ne\varnothing$. This requirement must be stated separately
when~${r=1}$, since facets are then zero-dimensional and are not covered by
the positive-dimensional face conditions below.

 \textbf{Hyperplane projections.}
The Zariski closure of an entire curve in $T_N\setminus D$ is a~translate of an algebraic subtorus contained in that complement.
Such subtori correspond to quotients of the character lattice $M$.
Combinatorially, such degenerations arise when the support $S$ collapses
under projection along rational directions.

For a~finite set $A \subset M_{\mathbb{R}}$ we define
\[
\mathcal{L}_A := \{a-b \mid a,b\in A\},\qquad
V_A := \operatorname{Span}_{\mathbb{R}}(\mathcal{L}_A),
\]
and let
\[
\mathcal{H}_A := \{ H \subset V_A \mid H \text{ is a codimension-one
	$\mathbb{R}$--subspace generated by elements of }\mathcal{L}_A \}.
\]
For $H\in\mathcal{H}_A$, we denote by
$
\phi_H\colon V_A \longrightarrow V_A/H \simeq \mathbb{R}
$
the canonical projection.
When $A=\tau\cap S$ for a~face $\tau$ of $P_D$, we write
$\mathcal{H}_{\tau\cap S}$ accordingly.
Throughout, a~face is allowed to be the whole polytope, so every facewise
hypothesis below includes $\tau=P_D$. For a~positive-dimensional face
$\tau$, we use the saturated character lattice~${M_\tau:=M\cap\operatorname{span}_{\mathbb R}(\tau-\tau)}$ of its torus
orbit; translating $\tau\cap S$ by any of its lattice points places the
support in $M_\tau$. Projection cardinalities are independent of this
translation.
For any nonempty finite support $A$ and $x\in A$, we use
$|\phi_H(A)|$ only as shorthand for $|\phi_H(A-x)|$.
This convention applies throughout, including the lemmas, tables, and
examples. Whenever an expression $|\phi_H(A-x)|$ is used without a~specified base point, $x$ denotes an arbitrary point of $A$.

 \textbf{Tiba's hyperbolicity criterion.}
With the above notation, Tiba's generic hyperbolicity theorem can be stated as
follows.

\begin{TibaTheorem}[Tiba {\cite{TibaYusaku13}}]
	\label{thm:Tiba13MainThm}
	Let $S\subset P_D\cap M$ be finite with $\mathfrak d_S\ne\varnothing$.
	Assume that for every positive-dimensional face $\tau$ of $P_D$:
	\begin{enumerate}
\itemsep=0pt
		\item[$(1)$] $\tau\cap S\neq\varnothing$, and the convex hull of $\tau\cap S$ has
		the same dimension as $\tau$;
		\item[$(2)$] for each hyperplane $H\in\mathcal{H}_{\tau\cap S}$ there exists
		$x\in\tau\cap S$ such that
		\[
		\bigl|\phi_H(\tau\cap S - x)\bigr| \ge \dim\tau+1.
		\]
	\end{enumerate}
	Then, for a~general divisor $D\in\mathfrak{d}_S$, the complement
	$T_N\setminus\Supp D$ is Kobayashi hyperbolically embedded in $X$.
\end{TibaTheorem}

\begin{Remark}
Theorem~\ref{thm:Tiba13MainThm} is the generic existence result used below.
We retain Tiba's incidence method and record its coefficient conditions in a~form suited to one-coefficient deformations and the calculations in
Section~\ref{sec:examples}.
\end{Remark}

\section{Effective proof of Tiba's hyperbolic embedding theorem}\label{sec:proof-main}

\subsection{Position of Tiba's theorem in this paper}\label{subsec:tiba-position}

We recall that the central generic hyperbolicity statement used throughout this paper
is due to Tiba and has already been stated as Theorem~\ref{thm:Tiba13MainThm}.
That theorem provides a~combinatorial criterion on a~finite support
$S\subset P_D\cap M$ ensuring that, for a~\emph{general} divisor
$D\in\mathfrak d_S$, the complement
$T_N\setminus\Supp D$ is Kobayashi hyperbolically embedded in the ambient
projective toric variety $X$.

The purpose of the present section is to give an effective and coefficient-wise
realization of Theorem~\ref{thm:Tiba13MainThm}.
More precisely, we reformulate Tiba's argument in a~way that
\begin{itemize}
\itemsep=0pt
	\item isolates an explicit constructible bad locus and a~proper Zariski
	closed exceptional locus~${\mathcal Y_S\subset\mathfrak d_S^*}$;
	\item characterizes membership in the constructible bad loci via
	translated subtori in torus strata;
	\item extracts a~standalone ``no entire curves'' criterion on $T_N$ that will be
	used later in the proof of Theorem~\ref{thm:2}.
\end{itemize}
No new generic hyperbolicity statement is claimed here; rather, this section
should be read as an effective proof and structural refinement of
Theorem~\ref{thm:Tiba13MainThm}.

\subsection{Setup and strategy}
\label{subsec:setup-strategy}

Fix $(X,T_N,D,P_D)$ as in Section~\ref{sec:toric-geometry}.
Let $S\subset P_D\cap M$ be a~finite support of full affine dimension with
$\mathfrak d_S\ne\varnothing$, and
use the parameter spaces $\mathfrak d_S\subset\mathfrak d_S^*$ defined in
Section~\ref{sec:toric-geometry}. The complete support~${S=P_D\cap M}$ will be used for the global locus $\mathcal Y\subset|D|$.

For a~nonzero coefficient vector $a$, the chosen toric trivialization on
$T_N$ represents $s_a$ by
\[
F_a(z):=\sum_{I\in S}a_Iz^I,\qquad z\in T_N.
\]
Its zero divisor is $D_T(a)$, and its closure is $D_{\mathrm{cl}}(a)$,
as defined in Section~\ref{sec:toric-geometry}. In the local
Laurent-polynomial arguments of this section, we abbreviate
$F=F_a$ and $D_0=D_T(a)$; this local notation does not denote a~divisor
on all of $X$. In the facewise argument, we use $D_X(a)$ explicitly.

Set
\[
r:=\dim_{\mathbb R} \mathrm{span}_{\mathbb R}(S-S)
=\dim \mathrm{Conv}(S).
\]
All subtori appearing in the argument have dimension at most $r$.

The proof proceeds in four steps:
\begin{enumerate}
\itemsep=0pt
	\item[(1)] reduction to one-dimensional subtori;
	\item[(2)] construction of explicit exceptional loci in coefficient space;
	\item[(3)] equivalence between bad coefficients and translated subtori;
	\item[(4)] exclusion of entire curves via the theorem of Dethloff--Lu.
\end{enumerate}

\subsection{Reduction to one-dimensional subtori}

\begin{Lemma}
	\label{lem:1d-reduction}
	Let $Z\subset T_N$ be any subset. Then $Z$ contains a~translate of a~positive-dimensional subtorus if and only if it contains a~translate of a~one-dimensional subtorus.
	In particular, the same equivalence holds for $Z=\Supp D_0$ and
	$Z=T_N\setminus \Supp D_0$.
\end{Lemma}

\begin{proof}
	If $Z$ contains a~translate of a~positive-dimensional subtorus, then it
	contains the translate of any one-dimensional subtorus inside it.
	Conversely, a~one-dimensional subtorus is itself positive-dimensional, so
	its translate already has the required form.
\end{proof}

\subsection{Exceptional loci in coefficient space}

Fix $S\subset P_D\cap M$ of full dimension $r$.
For each $H\in\mathcal{H}_S$, choose a~splitting
\[
M \cong (H\cap M)\oplus \big(M/(H\cap M)\bigr),
\qquad M/(H\cap M)\cong\mathbb{Z}.
\]
With respect to the induced coordinates
$z\leftrightarrow(u_1,\dots,u_{r-1},u_r)$ on $T_N$, we may rewrite
\[
F(z)=\sum_{i=1}^{\ell_H} Q_{H,i}(u_1,\dots,u_{r-1})u_r^{d_i},
\]
where $\ell_H:=|\phi_H(S)|$ and each $Q_{H,i}$ is linear in the coefficients
$(a_I)_{I\in S}$.

\begin{Definition}[bad and exceptional loci for $(S,H)$]\label{def:YSH}
	For $1\le j\le \ell_H$, define \smash{$B_{S,H,j}\subset\mathbb P^{|S|-1}$} to be the
	set of projective coefficient classes $[a]$, represented by
	$a=(a_I)_{I\in S}$, for which there exists~\smash{$c\in\bigl(\mathbb C^*\bigr)^{r-1}$} such that $
	Q_{H,k}(c)=0$ for all $k\neq j$.
	Set
	\[
	B_{S,H}:=\bigcup_{j=1}^{\ell_H} B_{S,H,j},
	\qquad
	B_S:=\bigcup_{H\in\mathcal{H}_S} B_{S,H}.
	\]
	We call $B_S$ the \emph{bad locus} and define the exceptional locus to be
	its reduced Zariski closure
	\[
	Y_S:=\bigl(\overline{B_S}^{ \mathrm{Zar}}\bigr)_{\mathrm{red}}
	\subset\mathbb P^{|S|-1}.
	\]
\end{Definition}

\begin{Lemma}
	\label{lem:YS-closed}
	Assume that $
	\ell_H=|\phi_H(S)|\ge r+1$ for every $H\in\mathcal H_S$.
	The bad locus $B_S$ is constructible, and $Y_S$ is a~proper Zariski closed
	reduced subscheme of
	$\mathbb P^{|S|-1}$ of codimension at least $1$.
\end{Lemma}

\begin{proof}
	Fix $H$ and $j$, and let
	\[
	Z_{S,H,j}:=
	\bigl\{(c,a)\in\bigl(\mathbb C^*\bigr)^{r-1}\times\mathbb P^{|S|-1}
	 \mid Q_{H,k}(c)=0\text{ for }k\ne j\bigr\}.
	\]
	This is closed in \smash{$\bigl(\mathbb C^*\bigr)^{r-1}\times\mathbb P^{|S|-1}$}.
	For a~fixed $c$, the equations $Q_{H,k}(c)=0$, $k\ne j$, are
	$\ell_H-1$ independent linear equations in the coefficients: their
	coefficient supports are nonempty and pairwise disjoint, and every monomial
	evaluated at $c\in\bigl(\mathbb C^*\bigr)^{r-1}$ is nonzero. Hence
	\[
	\dim Z_{S,H,j}\le (r-1)+(|S|-\ell_H)
	=|S|-\ell_H+r-1.
	\]
	The projection to coefficient space is not proper, so its image need not be
	closed. Chevalley's theorem does, however, show that
	$B_{S,H,j}=\operatorname{pr}_2(Z_{S,H,j})$ is constructible, and the same
	dimension estimate holds for its Zariski closure. Since
	$\ell_H\ge r+1$, that closure has dimension at most~${|S|-2}$, whereas
	$\dim\mathbb P^{|S|-1}=|S|-1$. Taking the finite union over $H$ and $j$
	proves the assertion for~\smash{$Y_S=\bigl(\overline{B_S}^{\mathrm{Zar}}\bigr)_{\mathrm{red}}$}.
\end{proof}

The translated-subtorus characterization in
Lemma~\ref{lem:YS-subtori-equiv} concerns $B_S$, not its closure $Y_S$.
We make no general assertion that $B_S=Y_S$, or that membership in $Y_S$
characterizes failure of hyperbolicity of the complement. A~translated
subtorus in the divisor itself is different from an entire curve in its
complement.

\subsection{Translated subtori and exceptional coefficients}

\begin{Lemma}
	\label{lem:YS-subtori-equiv}
	Let \smash{$(a_I)_{I\in S}\in\mathbb P^{|S|-1}$} and
	$D_0=\divv\bigl(\sum_{I\in S} a_I z^I\bigr)$.
	Then $(a_I)\in B_S$ if and only if either $\Supp D_0$ or
	$T_N\setminus\Supp D_0$ contains a~translate of a~one-dimensional subtorus
	of $T_N$.
\end{Lemma}

\begin{proof}
	Suppose first that $a\in B_{S,H,j}$, and choose $c$ as in
	Definition~\ref{def:YSH}. On the translated one-dimensional subtorus
	obtained by fixing $(u_1,\dots,u_{r-1})=c$, the restriction of $F$ is
	\[
	F(c,u_r)=Q_{H,j}(c)u_r^{d_j}.
	\]
	If $Q_{H,j}(c)=0$, this subtorus lies in $\Supp D_0$; otherwise, it is
	disjoint from $\Supp D_0$.

	Conversely, write a~translated one-dimensional subtorus as
	$c\lambda_v\bigl(\mathbb C^*\bigr)$, where $v\in N$ is primitive. Grouping the
	exponents according to $m=\langle I,v\rangle$ gives
	\[
	F(c\lambda_v(t))
	=
	\sum_m\biggl(\sum_{\langle I,v\rangle=m}a_Ic^I\biggr)t^m.
	\]
	If the translate is contained in $\Supp D_0$, every displayed coefficient
	vanishes. If it is contained in the complement, this Laurent polynomial
	has no zero in $\mathbb C^*$ and is therefore a~single monomial. Thus in
	either case at most one layer coefficient is nonzero.

	For completeness, the standard finite-hyperplane reduction is as follows.
	Put $S_m=\{I\in S\mid \allowbreak \langle I,v\rangle=m\}$ and
	$b_m=\sum_{I\in S_m}a_Ic^I$. Let $W\subset v^\perp$ be the span of
	the differences within the sets $S_m$. Since $S-S$ spans $M_{\mathbb R}$,
 a~basis of $W$ can be extended by differences from $S-S$ to a~hyperplane~${H\in\mathcal H_S}$. Each $H$-layer $L_d$ is a~union of sets $S_m$, so
	evaluation at the fixed base point $c$ gives
	\[
	\sum_{I\in L_d}a_Ic^I=\sum_{S_m\subset L_d}b_m.
	\]
	At most one of these sums is nonzero. In the coordinates adapted to $H$,
	write $c=(c_H,c_r)$. The left-hand side is
	$Q_{H,d}(c_H)c_r^d$, where $Q_{H,d}$ denotes $Q_{H,i}$ for the unique
	index $i$ satisfying $d_i=d$. Since $c_r\ne0$, at most one specialized
	$H$-layer polynomial is nonzero. Therefore,
	$a\in B_{S,H,j}\subset B_S$ for a~suitable $j$. This is the
	hyperplane-reduction argument used in~\cite[Lemma~2]{TibaYusaku13}.
\end{proof}

\subsection[The Dethloff--Lu theorem and a~''no entire curves'' criterion on T\_N]{The Dethloff--Lu theorem and a~``no entire curves'' criterion on $\boldsymbol{T_N}$}

\begin{Theorem}[Dethloff--Lu {\cite[Theorem~5.1]{DethloffLu2001}}]
	\label{thm:DethloffLu}
	Let $G$ be a~semi-abelian variety and let $D$ be a~reduced algebraic
	divisor on $G$. If $f\colon \mathbb C\to G$ is holomorphic, then the Zariski
	closure $X(f)$ of its image in $G$ is a~translate of an algebraic subgroup
	of $G$. If in addition $f(\mathbb C)\subset G\setminus D$, then~${X(f)\cap D=\varnothing}$.
\end{Theorem}

For an algebraic torus, the same degeneracy statement can also be obtained
from the classical Borel lemma. The semi-abelian formulation above applies
directly to the setting needed here.

The bad-locus construction and Lemma~\ref{lem:YS-subtori-equiv} depend
only on a~finite support of full affine dimension in the character lattice.
The following local statements therefore require neither a~nonempty global
coefficient chart nor conditions on proper faces. They apply also in the
affine lattice of any face after translating the support by a~lattice point.

\begin{Lemma}
	\label{lem:no-entire-curves-in-torus}
	Let $S\subset M$ be finite of full affine dimension $r$, and let
	$D_0=\divv\bigl(\sum_{I\in S}a_Iz^I\bigr)\subset T_N$.
	If $(a_I)\notin B_S$, then there exists no nonconstant holomorphic map
	from $\mathbb C$ to either $\Supp D_0$ or~${T_N\setminus\Supp D_0}$.
\end{Lemma}

\begin{proof}
	Suppose first that a~nonconstant $f\colon \mathbb C\to\Supp D_0$ exists.
	By Theorem~\ref{thm:DethloffLu}, the Zariski closure of $f(\mathbb C)$ is a~translate of a~positive-dimensional algebraic subtorus contained in
	$\Supp D_0$. If instead
	$f\colon \mathbb C\to T_N\setminus\Supp D_0$ is nonconstant, the same theorem
	shows that its Zariski closure is a~translate of a~positive-dimensional
	algebraic subtorus disjoint from $\Supp D_0$. In either case,
	Lemma~\ref{lem:1d-reduction} produces a~translated one-dimensional
	subtorus in the corresponding set, and
	Lemma~\ref{lem:YS-subtori-equiv} gives $(a_I)\in B_S$, a~contradiction.
\end{proof}

\begin{Corollary}[no entire curves in the torus complement]
	\label{thm:no-entire-curves-torus}
	Let $S\subset M$ be finite of full affine dimension $r$.
	Let $D_0=\divv\bigl(\sum_{I\in S} a_I z^I\bigr)\subset T_N$.
	If $(a_I)\notin B_S$, then $T_N\setminus\Supp D_0$ is Brody hyperbolic, i.e., there is no nonconstant holomorphic map $\C\to T_N\setminus\Supp D_0$.
\end{Corollary}

\begin{proof}
	This is the complement case of Lemma~\ref{lem:no-entire-curves-in-torus}.
\end{proof}

\subsection{Completion of the proof of Tiba's theorem}

\begin{Lemma}
\label{lem:facewise-hyperbolic-embedding}
Assume the hypotheses of Theorem~{\rm\ref{thm:Tiba13MainThm}}, and let
$D_X(a)\in\mathfrak d_S$. For every positive-dimensional face $\tau$,
put $S_\tau=\tau\cap S$ and let
$\rho_\tau(a)=(a_I)_{I\in S_\tau}$ be the restricted coefficient vector.
If, for every such face,
$
\rho_\tau(a)\ne0$, $
[\rho_\tau(a)]\notin B_{S_\tau}$,
then $T_N\setminus\operatorname{Supp}D_T(a)$ is Kobayashi
hyperbolically embedded in $X$.
\end{Lemma}

\begin{proof}
We want to verify Tiba's boundary criterion using
$D_{\mathrm{cl}}(a)$, the closure of the torus divisor.
Since $D_X(a)\in\mathfrak d_S$, this closure equals $D_X(a)$.

Fix a~positive-dimensional face $\tau$ and let $T_\tau\subset X$ be the
corresponding torus orbit. Choose $I_\tau\in S_\tau$. The monomial
section indexed by $I_\tau$ is nonvanishing on $T_\tau$ and trivializes
$\mathcal O_X(D)|_{T_\tau}$. In this trivialization the restricted
section is represented by
\[
F_{\tau,a}(z)=\sum_{I\in S_\tau}a_Iz^{I-I_\tau},
\qquad I-I_\tau\in M_\tau.
\]
The hypothesis $\rho_\tau(a)\ne0$ gives
$F_{\tau,a}\not\equiv0$.
Lemma~\ref{lem:no-entire-curves-in-torus}, applied in the character
lattice $M_\tau$, excludes nonconstant entire curves both in its zero
set and in its complement in $T_\tau$.
The full face $\tau=P_D$ treats the dense torus. Zero-dimensional orbits
are points, and therefore need no further condition. Hence all the
stratum conditions of~\cite[Lemma~3]{TibaYusaku13} hold for the reduced
support of $D_T(a)$ and its closure.
The conclusion of~\cite[Lemma~3]{TibaYusaku13} is the required
hyperbolic embedding.
\end{proof}

\begin{proof}[Proof of Theorem~\ref{thm:Tiba13MainThm}]
	For every positive-dimensional face $\tau$, put
	$S_\tau=\tau\cap S$ and let
	$\rho_\tau\colon\mathbb C^S\to\mathbb C^{S_\tau}$ be coordinate restriction.
	If $\widehat{Y}_{S_\tau}\subset\mathbb C^{S_\tau}$ denotes the affine cone
	over $Y_{S_\tau}$, including its vertex, define
	\[
	\mathcal Y_S:=
	\bigl(\mathfrak d_S^*\setminus\mathfrak d_S\bigr)
	\cup
	\bigcup_{\dim\tau>0}
	\mathbb P\bigl(\rho_\tau^{-1}\bigl(\widehat{Y}_{S_\tau}\bigr)\bigr)
	\subset\mathfrak d_S^*.
	\]
	Each term is closed and proper. Since only finitely many faces occur,
	$\mathcal Y_S$ is a~proper Zariski closed subset of
	$\mathfrak d_S^*$. For every $D_X(a)\notin\mathcal Y_S$, we have
	$D_X(a)\in\mathfrak d_S$. Moreover, every positive-dimensional face
	vector $\rho_\tau(a)$ is nonzero: otherwise it would lie at the vertex
	of the affine cone~$\widehat Y_{S_\tau}$ and the coefficient class
	would belong to~$\mathcal Y_S$. Its projective class lies outside
	$Y_{S_\tau}$ and hence outside~$B_{S_\tau}$.
	Lemma~\ref{lem:facewise-hyperbolic-embedding} therefore proves the
	generic assertion of Theorem~\ref{thm:Tiba13MainThm}. For the complete support~${S=P_D\cap M}$, we henceforth write $
	\mathcal Y:=\mathcal Y_{P_D\cap M}\subset |D|$.
\end{proof}

\section{Stability under one-parameter deformations}
\label{sec:deformation-stability}

In Section~\ref{sec:proof-main}, we gave an effective proof of
Tiba's hyperbolic embedding theorem (Theorem~\ref{thm:Tiba13MainThm}) by
constructing a~proper Zariski closed exceptional locus
$
\mathcal{Y}\subset |D|
$
in the coefficient parameter space.
As a~consequence, for every divisor
$D_0\in \mathfrak d_D\setminus\mathcal{Y}$,
the natural inclusion
$
T_N\setminus \Supp(D_0\cap T_N)\hookrightarrow X
$
is a~Kobayashi hyperbolic embedding.

The purpose of the present section is twofold.
First, we establish a~general deformation-stability statement:
along an algebraic one-parameter family whose image is not contained in
$\mathcal Y$, hyperbolic embedding holds outside a~finite set of parameters
(see Theorem~\ref{thm:1}).
Second, under a~stronger lattice hypothesis, we prove an exception-free
statement for a~distinguished class of affine one-parameter families, yielding
Theorem~\ref{thm:2}.
The latter requires a~refined analysis of borderline degenerations, which
motivates the introduction of an enlarged exceptional locus.

\subsection{Need for a~larger exceptional locus}

The constructible bad locus underlying $\mathcal{Y}$ detects degenerations in which, along some hyperplane
direction, the defining Laurent polynomial has \emph{at most one nonzero
layer} after specializing transverse variables. If no layer remains, the
corresponding translated subtorus lies in the divisor; if one remains, it
lies in the complement.
The exponential parametrization of any one-dimensional subtorus then gives
nonconstant entire curves. Theorem~\ref{thm:DethloffLu} is used in the
opposite direction, to pass from an entire curve to an algebraic-subgroup
translate.

For one-monomial deformations of the form
\[
D_t
=
\divv \biggl(
\sum_{I\in S\setminus\{I_0\}} a_I z^I + t z^{I_0}
\biggr),
\qquad t\in\mathbb C,
\] a~specialization with exactly two nonzero layers need not belong to the
ordinary bad incidence locus for that direction. However, changing one
coefficient may eliminate one of these two layers. To guarantee avoidance
of the ordinary bad locus throughout the deformation, we require at least
three nonzero layers initially. The enlarged locus below records the
failure of this stronger initial condition.

\subsection{Definition of the enlarged exceptional locus}

Fix $S\subset P_D\cap M$ satisfying the hypotheses of Theorem~\ref{thm:2}.
For each $H\in \mathcal{H}_S$, recall the decomposition
\[
\sum_{I\in S} a_I z^I
=
\sum_{i=1}^{\ell_H}
Q_{H,i}(u_1,\dots,u_{r-1}) u_r^{d_i},
\]
where $\ell_H = |\phi_H(S-x)|$ and each $Q_{H,i}$ is a~Laurent polynomial
in $(u_1,\dots,u_{r-1})$.

\begin{Definition}[enlarged bad and exceptional loci]
	\label{def:YSprime}
	For each $H \in \mathcal H_S$ and each pair of distinct indices
	$1 \le j_1 < j_2 \le \ell_H$, let
	\smash{$B'_{S,H,j_1,j_2} \subset \mathbb{P}^{|S|-1}$}
	be the set of projective coefficient classes $[a]$, represented by
	$a=(a_I)_{I\in S}$, such that there exists
	\smash{$(c_1,\dots,c_{r-1}) \in \bigl(\mathbb{C}^*\bigr)^{r-1}$} satisfying
	\[
	Q_{H,i}(c_1,\dots,c_{r-1}) = 0
	 \qquad \text{for all } \quad i \neq j_1,j_2.
	\]
	Define
	\[
	B'_{S,H}:=\bigcup_{j_1<j_2} B'_{S,H,j_1,j_2},
	\qquad
	B'_S := \bigcup_{H\in\mathcal H_S} B'_{S,H},
	\qquad
	Y'_S:=\bigl(\overline{B'_S}^{ \mathrm{Zar}}\bigr)_{\mathrm{red}}.
	\]
\end{Definition}

By construction, $(a_I)_{I\in S} \in B'_S$ if and only if,
for some hyperplane direction, the defining polynomial degenerates to one
with at most two nonzero layers after specialization of transverse
variables.
In particular,
$
B_S \subset B'_S$ and $ Y_S \subset Y'_S$.

\subsection[Algebraicity and codimension of the enlarged exceptional locus]{Algebraicity and codimension of $\boldsymbol{Y'_S}$}

\begin{Lemma}
	\label{lem:YSprime-codim}
	The set $B'_S$ is constructible, and
	$Y'_S \subset \mathbb{P}^{|S|-1}$ is a~proper Zariski closed reduced
	subscheme of codimension at least\/ $1$.
\end{Lemma}

\begin{proof}
	Fix $H \in \mathcal H_S$ and indices $j_1 < j_2$.
	Consider the closed subscheme
	\[
	Z'_{S,H,j_1,j_2}
	\subset \bigl(\mathbb{C}^*\bigr)^{r-1} \times \mathbb{P}^{|S|-1}
	\]
	defined by $
	Q_{H,i}(u_1,\dots,u_{r-1}) = 0$ for all $ i \neq j_1,j_2$.
	For fixed $(u_1,\dots,u_{r-1})$, the fiber is cut out by
	$\ell_H-2$ independent linear equations in $\mathbb{P}^{|S|-1}$,
	hence has dimension $
	|S|-1-(\ell_H-2)=|S|-\ell_H+1$.
	Since $\ell_H\ge r+2$ under the hypotheses of
	Theorem~\ref{thm:2}, we obtain
	\[
	\dim Z'_{S,H,j_1,j_2}
	\le (r-1)+(|S|-\ell_H+1)\le |S|-2.
	\]
	The projection to $\mathbb{P}^{|S|-1}$ is not proper, so we do not claim
	that its image is closed. Chevalley's theorem shows that the image
	\smash{$B'_{S,H,j_1,j_2}$} is constructible, while the displayed estimate gives
\smash{$\dim\overline{B'_{S,H,j_1,j_2}}^{\mathrm{Zar}}
	\le |S|-2$}.
	Taking the finite union over $H$, $j_1$, $j_2$ proves the assertion for its
	reduced Zariski closure~$Y'_S$.
\end{proof}

\begin{Definition}[facewise enlarged exceptional locus]
	\label{def:global-enlarged-locus}
	For every positive-dimensional face $\tau$ put~${S_\tau=\tau\cap S}$,
	let $\rho_\tau\colon\mathbb C^S\to\mathbb C^{S_\tau}$ be coordinate
	restriction, and let $\widehat{Y}'_{S_\tau}$ be the affine cone over~$Y'_{S_\tau}$. We set
	\[
	\mathcal Y'_S:=
	\bigl(\mathfrak d_S^*\setminus\mathfrak d_S\bigr)
	\cup
	\bigcup_{\dim\tau>0}
	\mathbb P\bigl(\rho_\tau^{-1}\bigl(\widehat{Y}'_{S_\tau}\bigr)\bigr)
	\subset\mathfrak d_S^*.
	\]
	Under the hypotheses of Theorem~\ref{thm:2}, this is a~proper Zariski
	closed subset by Lemma~\ref{lem:YSprime-codim}.
\end{Definition}

\subsection[Proof of Theorem 1.1]{Proof of Theorem~\ref{thm:1}}
\label{subsec:proof-finite-exceptions}

\begin{proof}[Proof of Theorem~\ref{thm:1}]
	Let $D_0\in \mathfrak d_D\setminus \mathcal Y$ and
	$\gamma\colon \mathbb{P}^1\to |D|$ be a~morphism with $\gamma(0)=D_0$.
	Set $D_t:=\gamma(t)$.
	
	The inverse image $\gamma^{-1}(\mathcal Y)$ is Zariski closed in
	$\mathbb P^1$. It is proper because $0\notin\gamma^{-1}(\mathcal Y)$,
	and every proper closed subset of $\mathbb P^1$ is finite. This argument
	also covers the case of a~constant family, when the inverse image is empty.
	
	For all but finitely many $t\in\mathbb{P}^1$ we therefore have
	$D_t\in \mathfrak d_D\setminus\mathcal Y$,
	and the hyperbolic embedding~$
	T_N\setminus \Supp(D_t\cap T_N)\hookrightarrow X$
	follows from Lemma~\ref{lem:facewise-hyperbolic-embedding}.
\end{proof}

 \textbf{Extension to algebraic curves.}
We prove that the exceptional parameter set is finite for an
irreducible algebraic curve base. Let $C$ be an irreducible complex
algebraic curve and $\gamma\colon C\to|D|$ a~morphism with
$\gamma(c_0)\notin\mathcal Y$ for some $c_0\in C$. Then
$
E_C:=\gamma^{-1}(\mathcal Y)\subsetneq C
$
is Zariski closed. Since $C$ is a~finite-type curve, $\dim E_C\le0$
and $E_C$ is finite.
For $c\in C\setminus E_C$,
Lemma~\ref{lem:facewise-hyperbolic-embedding} gives hyperbolic embedding. Neither completeness nor smoothness of $C$ is
needed. For a~reducible curve, apply this argument separately to each
irreducible component whose image is not contained in $\mathcal Y$.
No conclusion is asserted on a~component mapped entirely into
$\mathcal Y$.

\subsection[Proof of Theorem 1.2]{Proof of Theorem~\ref{thm:2}}

\begin{proof}
	Let $a=(a_I)_{I\in S}$ satisfy $[a]\notin\mathcal Y'_S$ and fix $I_0\in S$, as in the
	statement. Consider
	\[
	F_t(z)=tz^{I_0}+\sum_{I\ne I_0}a_Iz^I,
	\qquad t\in\mathbb C.
	\]
By Definition~\ref{def:global-enlarged-locus}, every original face
vector $\rho_\tau(a)$ is nonzero and its projective class avoids
$Y'_{S_\tau}$ and hence $B'_{S_\tau}$. Indeed, the affine cones in
that definition include their vertices.
	Fix a~positive-dimensional face $\tau$ and write $S_\tau=\tau\cap S$.
	If $I_0\notin S_\tau$, the face polynomial is independent of~$t$ and its
	coefficient vector lies outside $B'_{S_\tau}$, hence outside
	$B_{S_\tau}$. Suppose that $I_0\in S_\tau$. For any~${H\in\mathcal H_{S_\tau}}$ and any specialization of the transverse
	variables, the assumption
	$a|_{S_\tau}\notin B'_{S_\tau}$ says that at least three layer
	polynomials are nonzero at the original coefficient vector. Varying the
	coefficient of $z^{I_0}$ changes only the single $H$-layer containing
	$I_0$. At least two of the previously nonzero layers are therefore
	unchanged and remain nonzero for every $t\in\mathbb C$. In particular,
	the specialized polynomial never collapses to at most one nonzero layer,
	so the coefficient vector of~$F_t|_{S_\tau}$ lies outside
	$B_{S_\tau}$ for every $t$.

We first prove the conclusion of Theorem~\ref{thm:2} for $r\ge2$.
Every facet then has positive dimension. On every
positive-dimensional face the layer argument leaves at least two nonzero
layers after every specialization, so the face coefficient vector is
nonzero and its projective class avoids the ordinary bad locus.
In particular, every facet restriction is nonzero, which means that~${D_t=D_X(a(t))\in\mathfrak d_S}$ for all $t$.
Lemma~\ref{lem:facewise-hyperbolic-embedding} therefore gives the asserted
hyperbolic embedding for all $t\in\mathbb C$ when $r\ge2$.

We now prove the rank-one case of Theorem~\ref{thm:2} directly.
For $r=1$, the toric curve $X$ is~$\mathbb P^1$.
There are no transverse variables, and avoidance of $B'_S$ means that
$F_{a_{I_0}}$ has at least three nonzero terms. Changing one coefficient
leaves at least two nonzero terms in $F_t$ for every $t$.
After factoring out its lowest Laurent monomial, write
$F_t(z)=z^mP_t(z)$, where $P_t(0)\ne0$ and~${\deg P_t\ge1}$.
Thus the zero set $Z_t\subset\mathbb C^*$ is finite and nonempty.
The complement is
\[
T_N\setminus\operatorname{Supp}D_T(a(t))
=\mathbb P^1\setminus\bigl(\{0,\infty\}\cup Z_t\bigr).
\] A~sphere punctured at at least three distinct points is hyperbolically
embedded in $\mathbb P^1$; see, for example, \cite[Introduction, Paragraph~3]{Zaidenberg1986}.
If an endpoint coefficient is deleted at $t=0$, the global divisor
$D_X(a(0))$ can acquire a~torus-invariant point and leave
$\mathfrak d_S$. The inclusion~${T_N\setminus\operatorname{Supp}D_T(a(0))\hookrightarrow X}$
still depends only on the torus restriction and its closure, so the
punctured-sphere argument proves the conclusion at $t=0$.
\end{proof}

\begin{Remark}[positioning of the results]
	Under the facewise lattice hypothesis, Theorem~\ref{thm:1} establishes
	deformation stability along projective one-parameter families whose
	images are not contained in $\mathcal Y$: hyperbolic embedding holds
	outside a~finite set of parameters.
	Theorem~\ref{thm:2} exploits a~stronger combinatorial hypothesis to obtain
	an exception-free statement for a~distinguished affine one-parameter
	family.
	The enlarged bad locus $B'_S$ is tailored to isolate the borderline two-term
	degenerations relevant to this refinement, and
	Corollary~\ref{thm:no-entire-curves-torus} is the local ``no entire curves''
	criterion that makes the argument uniform along the affine family.
\end{Remark}

\section{Similarity deformations and degeneracy parameters}
\label{sec:similarity-exceptions}

In Theorem~\ref{thm:1}, we proved that for any one-parameter deformation
$\gamma\colon \mathbb{P}^1\to |D|$ with $\gamma(0)=D_0\in \mathfrak d_D\setminus\mathcal Y$,
the exceptional parameter set
$E_{\mathcal Y}(\gamma):=\{t\in\mathbb{P}^1\mid \gamma(t)\in \mathcal Y\}$
is finite. In general, $E_{\mathcal Y}(\gamma)$ need not be empty.
If $E_{\mathrm{fail}}(\gamma)$ denotes the parameters at which the dense
torus complement is not hyperbolically embedded in $X$, our construction
proves only
$
E_{\mathrm{fail}}(\gamma)\subseteq E_{\mathcal Y}(\gamma)$.
We do not identify $\mathcal Y$ with the locus of all nonhyperbolic
complements, nor assert equality of these two parameter sets.

For deformations arising from \emph{diagonal torus reparametrizations}
(``similarity''), invariance will show that this intrinsic exceptional set is
actually empty on the domain where the reparametrization is defined. The
quantity to be minimized is therefore the degeneracy set of the rational
scaling functions, and we relate its optimal size to a~simple invariant of
$\lambda\in\bigl(\mathbb{C}^*\bigr)^r$.

\subsection{Similarity via the diagonal torus action}

Let $\lambda=(\lambda_1,\dots,\lambda_r)\in\bigl(\mathbb{C}^*\bigr)^r$. The diagonal map
$
\Delta_\lambda\colon T_N\longrightarrow T_N$, $
(z_1,\dots,z_r)\longmapsto (\lambda_1 z_1,\dots,\allowbreak\lambda_r z_r)
$
is an automorphism of $T_N$.

\begin{Definition}[similar divisors and Laurent polynomials]
	\label{def:similar}
	Let $D_p,D_q\in\mathfrak d_D$ be effective divisors. We call them
	\emph{similar} if, for some $\lambda\in\bigl(\mathbb C^*\bigr)^r$, $
	D_q=\Delta_\lambda^*D_p$,
	where the torus action is extended to~$X$. In particular,
	$\Supp(D_q\cap T_N)=\Delta_{\lambda^{-1}}(\Supp(D_p\cap T_N))$.
	Using the fixed toric trivialization of~$\mathcal O_X(D)$ on $T_N$,
	the corresponding Laurent equations satisfy $
	q(z)=c p(\lambda z)$, $ c\in\mathbb C^*$.
	Indeed, the corresponding global sections have the same divisor and
	hence differ by a~nonzero constant. After rescaling $q$, we use the
	normalized expression $q(z)=p(\lambda z)$ throughout this section.
	Working with effective divisors retains multiplicities; for arbitrary
	equations of divisors on the torus alone, the ambiguity would instead
	be a~Laurent unit $c z^m$.
\end{Definition}

If \smash{$p(z)=\sum_{I\in S} a_I z^I$}, then similarity by $\lambda$ gives
\[
q(z)=p(\lambda z)=\sum_{I\in S} a_I \lambda^I z^I,
\qquad
\lambda^I:=\lambda_1^{i_1}\cdots\lambda_r^{i_r}.
\]

\subsection{Similarity preserves membership in the exceptional locus}

Recall that the exceptional locus $\mathcal Y\subset|D|$
constructed in the proof of Theorem~\ref{thm:Tiba13MainThm}
is defined by facewise conditions detecting translates of one-dimensional subtori
contained either in $\Supp D_0$ or in its complement.
More precisely, $\mathcal Y$ is the reduced Zariski closure of the facewise
bad loci described by Lemma~\ref{lem:YS-subtori-equiv}, together with the
boundary of $\mathfrak d_D$.

\begin{Lemma}[similarity invariance of $\mathcal Y$]
	\label{lem:similarity-invariance}
	Let $D_p,D_q\in\mathfrak d_D$ be divisors defined on $T_N$ by similar Laurent polynomials
	$p$ and $q=p(\lambda z)$.
Then
$D_p\in\mathcal Y \Longleftrightarrow D_q\in\mathcal Y$.
\end{Lemma}

\begin{proof}
 A~translate of a~one-dimensional subtorus in $T_N$ is carried by $\Delta_\lambda$
	to a~translate of a~one-dimensional subtorus, and the same holds after restricting to any
	torus stratum. Hence every facewise bad locus $B_{S_\tau}$ is invariant.
We verify that the torus action preserves every part of the global
exceptional locus. Since the action is an algebraic automorphism, it preserves
the reduced Zariski closures $Y_{S_\tau}$ and their inverse images under
coefficient restriction. It multiplies each coefficient $a_I$ by the
nonzero number $\lambda^I$. Therefore, a~restricted coefficient vector
vanishes before the action if and only if it vanishes after the action.
In particular, the coefficient-chart boundary
$\mathfrak d_D^*\setminus\mathfrak d_D$ is invariant. All pieces in the
construction of $\mathcal Y$ are consequently invariant, proving the
claim.
\end{proof}

\subsection[Canonical similarity paths and a ''degeneracy set'']{Canonical similarity paths and a~``degeneracy set''}

Assume $p$ and $q=p(\lambda z)$ are similar.
Consider the rational functions
$\Lambda_i(s:t):=\frac{s+\lambda_i t}{s+t}$, $ 1\le i\le r$,
on $\mathbb{P}^1$, and define (on their common domain of definition)
\[
p_{(s:t)}(z)
:=p\bigl(\Lambda_1(s:t)z_1,\dots,\Lambda_r(s:t)z_r\bigr).
\]
Let $U\subset \mathbb P^1$ be the maximal open subset where all $\Lambda_i$ are defined and
take values in $\mathbb C^*$ (so that the diagonal scaling is non-degenerate).
Then $p_{(s:t)}$ defines a~morphism $\gamma\colon U\to\mathfrak d_D$ by~${(s:t)\mapsto D_{(s:t)}:=\divv(p_{(s:t)})}$.

We have $
p_{(1:0)}(z)=p(z)$,
$
p_{(0:1)}(z)=p(\lambda z)=q(z)$,
so $\gamma$ connects $D_p$ and $D_q$ inside $U$.

On $U$, define the exceptional preimage by
$
E_{\mathcal Y}(\gamma):=\gamma^{-1}(\mathcal Y)$.
We reserve $E_{\mathrm{fail}}(\gamma)$ for the parameters in $U$ where
hyperbolic embedding actually fails. The degeneracy set of the chosen
scaling functions $\Lambda=(\Lambda_1,\dots,\Lambda_r)$ is
\[
E_{\mathrm{deg}}(\Lambda):=\mathbb P^1\setminus U
=\bigcup_{i=1}^r \Lambda_i^{-1}(\{0,\infty\}).
\]
For the above canonical choice, one checks
\[
E_{\mathrm{deg}}(\Lambda)\subseteq \{(1:-1)\}\cup\{(\lambda_i:-1)\mid 1\le i\le r\},
\qquad \text{hence} \quad |E_{\mathrm{deg}}(\Lambda)|\le r+1.
\]
Every $D_{(s:t)}$ with $(s:t)\in U$ is similar to $D_p$. Hence
Lemma~\ref{lem:similarity-invariance} gives the sharper dichotomy
\[
D_p\notin\mathcal Y \Longrightarrow
\gamma^{-1}(\mathcal Y)=\varnothing,
\qquad
D_p\in\mathcal Y \Longrightarrow
\gamma(U)\subset\mathcal Y.
\]
Thus, when $D_p\notin\mathcal Y$, the only unavailable parameters of this
rational presentation are the points of $E_{\mathrm{deg}}(\Lambda)$; there are
no intrinsic exceptional parameters on $U$.

\subsection{Similar deformations and a~minimality problem}

The previous construction is a~special case of the following.

\begin{Definition}[similar deformations]
	\label{def:similar-deformation}
	Fix $\lambda=(\lambda_1,\dots,\lambda_r)\in\bigl(\mathbb{C}^*\bigr)^r$.
 A~\emph{similar deformation} joining $p$ and $p(\lambda z)$ is given by nonzero rational functions
	$\Lambda_i\in\mathbb C\bigl(\mathbb P^1\bigr)^*$ $(1\le i\le r)$ such that
\begin{enumerate}\itemsep=0pt
		\item[(1)] $\Lambda_i(1:0)=1$ and $\Lambda_i(0:1)=\lambda_i$;
		\item[(2)] on the common open set $U\subset\mathbb P^1$ where all $\Lambda_i$ take values in $\mathbb C^*$,
		the family
		\[
		p_{(s:t)}(z)=p\bigl(\Lambda_1(s:t)z_1,\dots,\Lambda_r(s:t)z_r\bigr)
		\]
		induces a~morphism $\gamma\colon U\to\mathfrak d_D$ connecting $D_p$ and $D_{p(\lambda z)}$.
	\end{enumerate}
	We define the degeneracy set
	\[
	E_{\mathrm{deg}}(\Lambda):=\mathbb P^1\setminus U=\bigcup_{i=1}^r \Lambda_i^{-1}(\{0,\infty\}).
	\]
Zeros and poles are taken after cancellation of common factors in each
function $\Lambda_i$.
This notation refers to the chosen rational scaling presentation, not just
to its induced map into the coefficient space. The endpoint scaling
$\lambda$ is fixed in the minimality problem below.
\end{Definition}

The resulting minimality problem is algebraic: realize the fixed scaling
$\lambda$ by rational functions~$\Lambda_i$ while minimizing the cardinality
of the union of their zero and pole sets, namely~${|E_{\mathrm{deg}}(\Lambda)|}$.

\subsection{Multiplicative rank and the optimal degeneracy size}

\begin{Definition}[multiplicative rank]
	\label{def:mult-rank}
	For $\lambda=(\lambda_1,\dots,\lambda_r)\in\bigl(\mathbb{C}^*\bigr)^r$, define the
	\emph{multiplicative rank} $\operatorname{mrk}(\lambda)$ to be the smallest integer $k\ge 0$
	for which there exist $\mu_1,\dots,\mu_k\in\mathbb{C}^*$ and integers~${c_{ij}\in\mathbb{Z}}$
	such that
	\[
	\lambda_i=\prod_{j=1}^k \mu_j^{c_{ij}},
	\qquad 1\le i\le r.
	\]
	By convention, $\operatorname{mrk}(1,\dots,1)=0$.
\end{Definition}

This invariant is not the torsion-free rank of the subgroup generated by
the coordinates. For example, if $\zeta\ne1$ is a~root of unity, then
$\operatorname{mrk}(\zeta,1,\dots,1)=1$.

\begin{Proposition}[optimal degeneracy size for similarity deformations]\label{prop:min-exceptions}
	Let $\lambda\in\bigl(\mathbb{C}^*\bigr)^r$ with $\lambda\neq (1,\dots,1)$ and set $k:=\operatorname{mrk}(\lambda)$.
	Among all similar deformations $($Definition~{\rm\ref{def:similar-deformation})} joining $p$ and $p(\lambda z)$,
	the minimal possible value of $|E_{\mathrm{deg}}(\Lambda)|$ equals $k+1$.
\end{Proposition}

\begin{proof}
	\emph{Step $1$: construction with $k+1$ degeneracy points.}
	Choose $\mu_1,\dots,\mu_k\in\mathbb{C}^*$ and integers~$c_{ij}$ such that
	$\lambda_i=\prod_{j=1}^k \mu_j^{c_{ij}}$ for all $i$.
	Define
	\[
	\Lambda_i(s:t):=\prod_{j=1}^k \left(\frac{s+\mu_j t}{s+t}\right)^{c_{ij}}.
	\]
	Then $\Lambda_i(1:0)=1$ and $\Lambda_i(0:1)=\lambda_i$.
	All zeros/poles of $\Lambda_i$ are contained in the finite set
	\[
	\{(1:-1)\}\cup\{(\mu_j:-1)\mid 1\le j\le k\},
	\]
	so $|E_{\mathrm{deg}}(\Lambda)|\le k+1$ for this deformation.
	
	\emph{Step $2$: lower bound $|E_{\mathrm{deg}}(\Lambda)|\ge k+1$.}
	Let a~similar deformation be given, and set
	\[
	E:=E_{\mathrm{deg}}(\Lambda)=\bigcup_{i=1}^r \Lambda_i^{-1}(\{0,\infty\}).
	\]
	Since $\lambda\ne(1,\dots,1)$, some $\Lambda_i$ is nonconstant and
	$E\ne\varnothing$. Choose $p_0\in E$ and a~coordinate $u$ on
	$\mathbb P^1\setminus\{p_0\}\cong\mathbb A^1$. The endpoints
	$A=(1:0)$ and $B=(0:1)$ lie outside $E$; put $a=u(A)$ and $b=u(B)$.
	After cancelling common factors, the zero and pole orders
	$m_{i,\alpha}\in\mathbb Z$ and the normalization $\Lambda_i(A)=1$ give
	\[
	\Lambda_i(u)=\prod_{\alpha\in u(E\setminus\{p_0\})}
	\left(\frac{u-\alpha}{a-\alpha}\right)^{m_{i,\alpha}}.
	\]
	Consequently,
	\[
	\lambda_i=\prod_{\alpha\in u(E\setminus\{p_0\})}
	\left(\frac{b-\alpha}{a-\alpha}\right)^{m_{i,\alpha}}.
	\]
	All displayed ratios are nonzero, and there are at most $|E|-1$ of them.
	Thus $\operatorname{mrk}(\lambda)\le |E|-1$.
	Therefore, $|E|\ge k+1$.
	
	Combining the two steps proves the optimality.
\end{proof}

\subsection[Consequence for similarity paths in Theorem 1.1]{Consequence for similarity paths in Theorem~\ref{thm:1}}

Let $p$ and $q$ be similar with $q(z)=p(\lambda z)$. By Lemma~\ref{lem:similarity-invariance},
either both divisors lie in $\mathcal Y$ or neither does. If $D_p\notin\mathcal Y$, then
for any (non-degenerate) similarity path $\gamma\colon U\to\mathfrak d_D$ connecting $D_p$ and~$D_q$,
Lemma~\ref{lem:similarity-invariance} implies
$\gamma^{-1}(\mathcal Y)=\varnothing$; Theorem~\ref{thm:1} is not needed for this conclusion.
For a~fixed endpoint scaling $\lambda\ne(1,\dots,1)$,
Proposition~\ref{prop:min-exceptions} shows that the minimum is
\[
\min |E_{\mathrm{deg}}(\Lambda)|=\operatorname{mrk}(\lambda)+1.
\]
The minimum is taken over rational scaling presentations with these
prescribed endpoint scalings. For $\lambda=(1,\dots,1)$, the constant
scalings $\Lambda_i\equiv1$ instead give minimum $0$.
	
% ======================================================
\section{Examples}
\label{sec:examples}

This section illustrates how the exceptional locus (globally $\mathcal Y\subset|D|$, and
in restricted supports~${\mathcal Y_S\subset\mathfrak d_S^*}$) and the deformation principle in Theorem~\ref{thm:1}
behave in concrete toric settings.
We focus on projective spaces (where coefficients admit explicit descriptions) and
Hirzebruch surfaces (where birational modifications reveal toric boundary phenomena).

\subsection{Projective spaces}
\label{ex_projspace}
\label{subsec:examples-projective}

\subsubsection[The case X=P\^{}2 and D=2H]{The case $\boldsymbol{X=\mathbb{P}^2}$ and $\boldsymbol{D=2H}$}
\label{subsubsec:P2-2H}

Let $X=\mathbb{P}^2$ with homogeneous coordinates $[x:y:z]$ and $D=2H$.
Then $P_D$ is the standard degree-$2$ simplex in $M_{\mathbb R}\cong\mathbb R^2$
\[
P_D=\bigl\{(a,b)\in\mathbb R^2\mid a\ge 0,\, b\ge 0,\, a+b\le 2\bigr\},
\]
and $|P_D\cap M|=6$. A~general divisor $D_0\sim 2H$ is a~plane conic
\[
\Supp D_0=\bigl\{a_{1} x^2+ a_{2} y^2 + a_{3} z^2+b_{1} yz + b_{2} xz + b_{3} xy = 0\bigr\},
\]
so $\mathfrak d_D^*=|D|\simeq\mathbb{P}^5$ with homogeneous coordinates
$(a_1,a_2,a_3,b_1,b_2,b_3)$.

The bad locus is closed and explicit in this case, so
$B_{P_D\cap M}=Y_{P_D\cap M}\subset\mathbb{P}^5$ can be written as a~union
of six hypersurfaces
\begin{align*}
B_{P_D\cap M}=Y_{P_D\cap M}
	={}& \{a_{1}=0\} \cup \{a_{2}=0\} \cup \{a_{3}=0\}\cup \{a_{1}b_{1}^2 + a_{2}b_{2}^2= b_1 b_2 b_3\}\\
	&\cup \{a_2 b_2 ^2 + a_3 b_3^2 = b_1 b_2 b_3\}\cup \{a_3 b_3^2 + a_1 b_1^2 = b_1 b_2 b_3\}.
\end{align*}

To verify this description, identify the six lattice points with
$(2,0)$, $(0,2)$, $(0,0)$, $(0,1)$, $(1,0)$, $(1,1)$. Up to sign, the six directions in
$\mathcal H_{P_D\cap M}$ are generated by
$
(1,0)$, $(0,1)$, $(1,-1)$, $(1,1)$, $(1,-2)$, $(2,-1)$.
The first three directions have three layers. For example, for the direction
generated by $(1,0)$, after setting $z=1$ the three layer polynomials are
\[
Q_0 (x)=a_1x^2+b_2x+a_3,\qquad
Q_1 (x)=b_3x+b_1,\qquad Q_2=a_2.
\]
The cases in which $Q_2$ is one of the vanishing layers lie in the coordinate
hyperplane $\{a_2=0\}$. A~common zero of $Q_0$ and $Q_1$ in
$\mathbb C^*$ implies
$
a_1b_1^2+a_3b_3^2-b_1b_2b_3=0$.
The other two three-layer directions give the remaining two coordinate
components and the remaining two cubic equations. Each of the last three
directions has five layers, with four singleton layers; its bad incidence
locus is contained in the coordinate components already listed. Thus the
bad locus is contained in the displayed closed union.

For the reverse inclusion, suppose first that $a_2=0$, so that
$p(x,y,1)=Q_0 (x)+yQ_1 (x)$. If either $Q_j$ is identically zero or has a~zero in $\mathbb C^*$, fixing a~suitable $x$ leaves at most one nonzero
layer. Otherwise both $Q_j$ are nonzero monomials. Then $p$ has two
monomials, and its restriction to a~one-dimensional subtorus in the kernel
of their character ratio is zero or a~single monomial. In all cases the
coefficient vector belongs to the bad locus. Permuting the variables
handles the other coordinate hyperplanes.
Outside these hyperplanes, consider
$a_1b_1^2+a_2b_2^2=b_1b_2b_3$. If~${b_1b_2\ne0}$, restricting to
$b_1y+b_2x=0$ and $z=1$ gives $p=a_3\ne0$. If $b_1b_2=0$, the equality
forces $b_1=b_2=0$; choose a~nonzero root $\rho$ of
$a_1\rho^2+b_3\rho+a_2=0$. Then $p(\rho u,u,1)=a_3$ for~${u\in\mathbb C^*}$. These are translated subtori in the complement.
Permuting variables handles the other two cubic equations. This proves
the displayed equality with the reduced closed union, including all
degenerate coefficient vectors.

Moreover, for any one-dimensional face $\tau$ of $P_D$, the set $\tau\cap M$ has
three lattice points and $(a,b,c)\in B_{\tau\cap M}$ if and only if at least two among
$a$, $b$, $c$ vanish. In particular, failure of an edge condition forces membership in $Y_{P_D\cap M}$,
so the global condition already controls all positive-dimensional faces.
We identify the global conic exceptional locus by checking the edge
contributions and the coefficient-chart boundary. An edge bad locus
requires at least two of its three coefficients to vanish. At least one of these is an endpoint coefficient
among $a_1$, $a_2$, $a_3$. The coefficient-chart boundary requires all three
coefficients on some edge to vanish and is contained in the same union of
vertex-coefficient hyperplanes. These hyperplanes already occur in
$Y_{P_D\cap M}$. Since the full face also occurs in the global
construction, we obtain
\smash{$
\mathcal Y=Y_{P_D\cap M}=B_{P_D\cap M}$} in the conic case.

\begin{Proposition}[a concrete sufficient condition for hyperbolic embedding]\label{prop:P2-2H-sufficient}
Let $D_0=\divv\bigl(a_{1} x^2+ a_{2} y^2 + a_{3} z^2+b_{1} yz + b_{2} xz + b_{3} xy\bigr)$.
	If $
	a_1\neq 0$, $ a_2\neq 0$, $ a_3\neq 0$,
	and
	\begin{align*}
		a_{1}b_{1}^2 + a_{2}b_{2}^2&\neq b_1 b_2 b_3,\qquad
		a_2 b_2 ^2 + a_3 b_3^2 \neq b_1 b_2 b_3,\qquad
		a_3 b_3^2 + a_1 b_1^2 \neq b_1 b_2 b_3,
	\end{align*}
	then $T_N\setminus \Supp D_0 \hookrightarrow\mathbb{P}^2$ is a~Kobayashi hyperbolic embedding.
\end{Proposition}

\begin{proof}
	The inequalities in the statement say exactly that the full coefficient
	vector lies outside~${B_{P_D\cap M}}$. On each edge of $P_D$, two endpoint
	coefficients are among $a_1$, $a_2$, $a_3$ and are nonzero. Hence at least two of
	the three edge coefficients are nonzero, so the restricted coefficient
	vector lies outside the one-dimensional face bad locus as well.
	Lemma~\ref{lem:facewise-hyperbolic-embedding} now gives the asserted
	hyperbolic embedding.
\end{proof}

Failure of one of these sufficient inequalities can genuinely produce a~nonhyperbolic complement. Two typical degenerations are as follows:
\begin{itemize}\itemsep=0pt
\item[$(i)$] \emph{Vanishing of a~diagonal coefficient.}
Assume $a_1=0$ and $b_2b_3\ne0$, and put
\[
q(y,z)=a_2y^2+b_1yz+a_3z^2.
\]
If $q(-b_2,b_3)\ne0$, then on the torus locus
$b_3y+b_2z=0$ one has
$
p(x,y,z)=q(y,z)\ne0
$
for every~${x\in\mathbb C^*}$. Thus the complement contains the translated
one-dimensional subtorus obtained by fixing $(y,z)=(-b_2,b_3)$ and varying
$x$, and consequently it is not Brody hyperbolic.

\item[$(ii)$] \emph{Satisfying a~``cubic'' equality.}
Assume $a_1,a_2,a_3\neq0$ and
$a_1b_1^2+a_2b_2^2=b_1b_2b_3$. If $b_1b_2\ne0$, then
\[
p(x,y,z)= (b_1y+b_2x)
\left(\frac{a_1}{b_2}x+\frac{a_2}{b_1}y+z\right)+a_3 z^2.
\]
On $b_1y+b_2x=0$ with $z=1$, the value of $p$ is the nonzero constant
$a_3$, so the complement again contains a~translated subtorus. The equality
also permits the previously omitted branch $b_1=b_2=0$. In that case choose
$\rho\in\mathbb C^*$ satisfying
$a_1\rho^2+b_3\rho+a_2=0$; then~${u\mapsto(\rho u,u,1)}$ is a~translated subtorus on which $p=a_3\ne0$.
\end{itemize}

Here is a~one-parameter deformation with finitely many exceptions. Fix
$b_1,b_2\in\mathbb C^*$ and consider
$
p_t(x,y,z)=x^2+y^2+z^2+b_1yz+b_2xz+t xy$.
We determine exactly the parameters at which hyperbolic embedding
fails. Define
\begin{gather}
g_1(t):=b_1^2+b_2^2-b_1b_2t,\qquad
g_2(t):=b_2^2+t^2-b_1b_2t,\qquad
g_3(t):=t^2+b_1^2-b_1b_2t.\label{eq:conic-failure-polynomials}
\end{gather}
The exact failure set is
\begin{equation}\label{eq:conic-failure-set}
E_{\mathrm{fail}}
=\{\infty\}\cup\{t\in\mathbb C\mid g_1(t)g_2(t)g_3(t)=0\}.
\end{equation}

We first prove hyperbolic embedding outside the set in
\eqref{eq:conic-failure-set}. If $g_1(t)g_2(t)g_3(t)\ne0$, all the
hypotheses of Proposition~\ref{prop:P2-2H-sufficient} hold with
$a_1=a_2=a_3=1$ and $b_3=t$, giving the required embedding.

We prove failure at each finite parameter in
\eqref{eq:conic-failure-set} by constructing a~translated subtorus.
For~${g_i(t)=0}$, use the corresponding map
\begin{gather}
\iota_1(u)=[b_1u:-b_2u:1],\qquad
\iota_2(u)=[1:b_2u:-tu],\nonumber\\
\iota_3(u)=[b_1u:1:-tu],\qquad u\in\mathbb C^*.\label{eq:conic-subtori}
\end{gather}
Using~\eqref{eq:conic-failure-polynomials}, we compute the polynomial
evaluations at the representatives in~\eqref{eq:conic-subtori}:
$p_t(\iota_i(u))=1+g_i(t)u^2=1$ if $g_i(t)=0$.
The first map has nonzero coordinates because $b_1b_2\ne0$.
For the other maps, $g_2(0)=b_2^2\ne0$ and $g_3(0)=b_1^2\ne0$
ensure $t\ne0$ on the corresponding zero sets.
Thus each relevant $\iota_i$ lies in the torus complement and
$s\mapsto\iota_i\bigl(e^s\bigr)$ is a~nonconstant entire curve there.

We finally treat $t=\infty$ in~\eqref{eq:conic-failure-set}.
The section at infinity is $xy$, which has no zero on $T_N$.
The torus complement is therefore $T_N$ itself and is not hyperbolic.
The three cases establish~\eqref{eq:conic-failure-set}.

\subsubsection[The case X=P\^{}2 and D=3H]{The case $\boldsymbol{X=\mathbb{P}^2}$ and $\boldsymbol{D=3H}$}
\label{subsubsec:P2-3H}

Now let $D=3H$. Then
\[
P_D=\{(m_1,m_2)\in M_{\mathbb R}\mid m_1\ge 0,\, m_2\ge 0,\, m_1+m_2\le 3\}
\]
and $|P_D\cap M|=10$. A~general cubic divisor $D_0\sim 3H$ is defined by
\begin{align*}
	p(x,y,z)
	={}& a_{03}y^3
	+ a_{02}y^2z + a_{12}xy^2+ a_{01}yz^2 + a_{11}xyz + a_{21}x^2y\\
	&+ a_{00}z^3 + a_{10}z^2x + a_{20}zx^2 + a_{30}x^3.
\end{align*}

The calculation is guided by the following observation. The arrangement
$\mathcal H_{P_D\cap M}$ contains finitely many linear directions $H\subset M_{\mathbb R}$,
and each $H$ yields a~determinantal-type bad condition~$B_H$ in coefficient space.
Rather than listing all components, we record usable sufficient conditions by grouping the
relevant $H$ into geometric families.

We first consider the facet-parallel lines. Define
\begin{gather*}
	H_x=\{(m_1,m_2)\mid m_1=0\},\qquad
	H_y=\{(m_1,m_2)\mid m_2=0\},\\
	H_z=\{(m_1,m_2)\mid m_1+m_2=0\}.
\end{gather*}
These equations specify vector directions through the origin; affine translates
of these lines cut out the corresponding layers in $P_D$.
The condition ``\smash{$(a_{ij})\notin B_{H_x}$}'' can be phrased as: for every $(y,z)\in\bigl(\mathbb C^*\bigr)^2$,
at least two of the following four expressions are nonzero
\begin{gather*}
	F_0(y,z):=a_{00}z^3+a_{01}yz^2+a_{02}y^2z+a_{03}y^3,\qquad
	F_1(y,z):=a_{10}z^2+a_{11}yz+a_{12}y^2,\\
	F_2(y,z):=a_{20}z+a_{21}y,\qquad
	F_3:=a_{30}.
\end{gather*}
Similarly, one obtains the corresponding lists for $H_y$ and $H_z$ by symmetry.

Three additional difference directions are represented by
\[
H_{xy}=\{m_1-m_2=0\},\qquad
H_{yz}=\{m_1+2m_2=0\},\qquad
H_{zx}=\{2m_1+m_2=0\}.
\]

\begin{Proposition}[a workable sufficient criterion for hyperbolic embedding in the cubic case]
	\label{prop:P2-3H-sufficient}
	Let \smash{$D_0=\divv\big(\sum_{i+j\le 3} a_{ij}x^i y^j z^{3-i-j}\bigr)$ in $\mathbb{P}^2$}.
	Assume the following ``non-concentration'' conditions:
	\begin{enumerate}\itemsep=0pt
		\item[\textup{(1)}] $($\textup{$H_x$}$)$ For every $(y,z)\in\bigl(\mathbb C^*\bigr)^2$, at least two among
		$F_0(y,z),F_1(y,z),F_2(y,z),F_3$ are nonzero.
		\item[\textup{(2)}] $($\textup{$H_y$}$)$ For every $(x,z)\in\bigl(\mathbb C^*\bigr)^2$, at least two among
		\[
\begin{aligned}
&a_{00}z^3+a_{10}z^2x+a_{20}zx^2+a_{30}x^3,\\
&a_{01}z^2+a_{11}xz+a_{21}x^2,\qquad a_{02}z+a_{12}x,\qquad a_{03}.
\end{aligned}
\]
		are nonzero.
		\item[\textup{(3)}] $($\textup{$H_z$}$)$ For every $(x,y)\in\bigl(\mathbb C^*\bigr)^2$, at least two among
		\[
\begin{aligned}
&a_{30}x^3+a_{21}x^2y+a_{12}xy^2+a_{03}y^3,\\
&a_{20}x^2+a_{11}xy+a_{02}y^2,\qquad a_{10}x+a_{01}y,\qquad a_{00}.
\end{aligned}
\]
		are nonzero.
		\item[\textup{(4)}] $($\textup{$H_{xy}$}$)$ For every $(x,y,z)\in\bigl(\mathbb C^*\bigr)^3$, at least two among
		\[
		a_{30}, a_{20}, a_{02}, a_{03},
		a_{10} z^2 + a_{21} xy, a_{00} z^2 + a_{11}xy,a_{01} z^2 + a_{12} xy
		\]
		are nonzero.
		\item[\textup{(5)}] $($\textup{$H_{yz}$}$)$ For every $(x,y,z)\in\bigl(\mathbb C^*\bigr)^3$, at least two among
		\[
		a_{03}, a_{12}, a_{10}, a_{00},
		a_{21} x^2 + a_{02} yz, a_{30} x^2 + a_{11}yz,a_{20} x^2 + a_{01} yz
		\]
		are nonzero.
		\item[\textup{(6)}] $($\textup{$H_{zx}$}$)$ For every $(x,y,z)\in\bigl(\mathbb C^*\bigr)^3$, at least two among
		\[
		a_{00}, a_{01}, a_{21}, a_{30},
		a_{02} y^2 + a_{10} zx, a_{03} y^2 + a_{11}zx,a_{12} y^2 + a_{20} zx
		\]
		are nonzero.
		\item[\textup{(7)}] $($\textup{face conditions}$)$ In each row of
		\[
		\begin{pmatrix}
			a_{03} & a_{12} & a_{21} & a_{30} \\
			a_{03} & a_{02} & a_{01} & a_{00} \\
			a_{30} & a_{20} & a_{10} & a_{00}
		\end{pmatrix}
		\]
		at least two entries are nonzero.
	\end{enumerate}
	Then $T_N\setminus \Supp D_0\hookrightarrow\mathbb{P}^2$ is a~Kobayashi hyperbolic embedding.
\end{Proposition}

\begin{proof}
	Every element of $\mathcal H_{P_D\cap M}$ is the span of a~primitive
	difference vector. Normalize such a~vector $v=(p,q)$ up to sign and count
	the levels of the quotient functional $
	\nu_v(i,j)=qi-pj$.
	Enumerating the primitive directions in
	$(P_D\cap M)-(P_D\cap M)$ gives twelve directions. They split as follows:
	\[
	\begin{array}{@{}ccc@{}}
\toprule
|\phi_H(P_D\cap M)|&\text{primitive generators of }H&
	\text{conditions used}\ \\ \midrule
	4&(1,0),\ (0,1),\ (1,-1)&\textup{(1)--(3)}\\
	7&(1,1),\ (1,-2),\ (2,-1)&\textup{(4)--(6)}\\
	9&(1,-3),\ (1,2),\ (2,-3),\ (2,1),\ (3,-2),\ (3,-1)&\textup{(7)}\\
\bottomrule
\end{array}
	\]
	For the first six directions, the displayed expressions in
	\textup{(1)--(6)} are precisely their layer polynomials. Thus the stated
	non-concentration conditions say that at least two specialized layers are
	nonzero.

	For each of the six nine-layer directions, the ten lattice points are
	distributed among eight singleton layers and one two-point layer. Condition
	\textup{(7)} forces at least three boundary coefficients to be nonzero:
	indeed, two lattice points cannot supply two nonzero entries on each of the
	three edges. If at least four coefficients in total are nonzero, at least
	two nonzero singleton layers remain, regardless of cancellation in the
	two-point layer. If exactly three coefficients in total are nonzero,
	condition \textup{(7)} requires these three points to supply at least six
	incidences with the edges. They must therefore be the three vertices.
	The vertices lie in distinct layers because none of these six directions
	is parallel to an edge, and there are no other nonzero coefficients to
	cancel them. Thus every specialization has at least two nonzero layers
	in each of the twelve directions. The same edge
	condition excludes the bad loci of the three one-dimensional faces.
	Lemma~\ref{lem:facewise-hyperbolic-embedding} proves the conclusion.
\end{proof}

\begin{Proposition}[uniform hyperbolicity along one-monomial deformations]
	\label{prop:P2-3H-uniform}
	Fix a~coefficient vector $(a_{ij})$ satisfying the conditions of
	Proposition~{\rm\ref{prop:P2-3H-sufficient}} with ``at least two'' replaced by
	``at least three'' in \textup{(1)}--\textup{(7)}. Then, for any lattice
	point $I=(i_0,j_0)\in P_D\cap M$, let
	$m_I=x^{i_0}y^{j_0}z^{3-i_0-j_0}$. Replacing its coefficient by $t$ gives
	the family
	\[
	D_t=\divv \biggl(t m_I+\sum_{\substack{i+j\le 3\\ (i,j)\neq I}} a_{ij}x^i y^j z^{3-i-j}\biggr).
	\]
	Then $T_N\setminus \Supp D_t\hookrightarrow\mathbb{P}^2$ is a~Kobayashi
	hyperbolic embedding for all $t\in\mathbb C$.
\end{Proposition}

\begin{proof}
	The proof of Proposition~\ref{prop:P2-3H-sufficient} already identifies all
	twelve full-face directions. Conditions \textup{(1)--(6)} give at least
	three nonzero layers for the six directions with four or seven layers. For
 a~nine-layer direction, condition \textup{(7)} with the strengthened bound
	forces at least six boundary coefficients to be nonzero. Indeed, if $v$ is
	the number of nonzero vertex coefficients and $e$ the number of nonzero
	coefficients in the relative interiors of the edges, then counting
	incidences with the three edges gives $2v+e\ge9$; hence
	$v+e\ge9-v\ge6$. Since there is only one two-point layer, at least four
	nonzero singleton layers remain, regardless of cancellation in that layer.
	In particular, at least three specialized layers are nonzero.
	The strengthened face conditions give
	the same three-layer bound on every edge.

	Varying the coefficient indexed by $I$ changes only one layer for each
	direction and for each face. At least two of the previously nonzero layers
	therefore remain unchanged for every $t\in\mathbb C$. Thus no member of the
	family enters an ordinary facewise bad locus, and
	Lemma~\ref{lem:facewise-hyperbolic-embedding} yields the conclusion for all~$t$.
\end{proof}

\begin{Proposition}[deformation to a~five-term cubic]
	\label{prop:P2-3H-five-term}
	There is a~nonempty Zariski open subset of the space of plane cubics such
	that every divisor $D_0$ in this subset admits a~sequence of deformations
	preserving hyperbolic embedding to a~divisor
	$D'$ defined by a~five-term cubic
	\[
	a_{30}x^3+a_{03}y^3+a_{10}xz^2+a_{02}y^2z+a_{00}z^3.
	\]
\end{Proposition}

\begin{proof}
	Starting with the ten lattice points of the degree-$3$ simplex, eliminate
	the coefficients in the order $
	a_{11}$, $ a_{01}$, $ a_{12}$, $ a_{20}$, $ a_{21}$.
	Equivalently, remove the lattice points
	$(1,1)$, $(0,1)$, $(1,2)$, $(2,0)$, $(2,1)$ in that order. The layer counts used below
	are obtained by the following finite enumeration. For each current support
	$S_k$, take every nonzero difference $v=(p,q)\in S_k-S_k$, replace it by
	its primitive representative up to sign, and count the distinct integers
	$qi-pj$ for $(i,j)\in S_k$. For an edge, the only quotient is the identity
	on that one-dimensional affine lattice, so its layer count is simply the
	number of remaining edge points. This gives:
	\[
	\begin{array}{@{}cccc@{}}
\toprule
\text{coefficient to be removed}&|S|\text{ before deletion}&
	\min_{H\in\mathcal H_S}|\phi_H(S)|&
	\min_{\text{edges }\tau}|S\cap\tau|\\ \midrule
	a_{11}&10&4&4\\
	a_{01}&9&4&4\\
	a_{12}&8&4&3\\
	a_{20}&7&4&3\\
	a_{21}&6&4&3\\
\bottomrule
\end{array}
	\]
	For example, for the initial ten-point support the twelve normalized
	directions are those listed in the proof of
	Proposition~\ref{prop:P2-3H-sufficient}; the three edge-parallel directions
	realize the minimum $4$. Repeating the same normalization after each
	listed deletion gives the five rows of the table, so no unlisted direction
	is being omitted.
	These are exactly the enlarged bounds $\dim\tau+2$ for the two-dimensional
	face and its one-dimensional faces. More explicitly, let $S_0$ be the
	ten-point support and let~$S_k$ $(0\le k\le4)$ be the support just before
	the $(k+1)$-st deletion. The coordinate restrictions~${\rho_k\colon \mathbb C^{S_0}\to\mathbb C^{S_k}}$ are surjective. Each affine
	cone \smash{$\widehat{\mathcal Y'_{S_k}}$} is proper and closed, so the set
	\[
	\mathbb C^{S_0}\setminus
	\left(\bigcup_{k=0}^{4}\rho_k^{-1}
	\bigl(\widehat{\mathcal Y'_{S_k}}\bigr) \cup
	\bigcup_{I\in S_0}\{a_I=0\}\right)
	\]
	is nonempty and Zariski open. Choose the initial coefficients once in
	this set. At each stage the coefficients still present are exactly
	$\rho_k(a)$: earlier deletions changed only coordinates outside $S_k$.
	Thus no new generic choice is made during the process, and
	Theorem~\ref{thm:2} applies successively. Each affine coefficient line
	preserves hyperbolic embedding from the current coefficient to zero.
	The remaining support is $
	\{(3,0),(0,3),(1,0),(0,2),(0,0)\}$,
	which gives the displayed five-term cubic. Its minimum full-face layer
	count is $3$, and its three edge counts are $3$, $3$, $2$, so it still satisfies
	the ordinary bounds of Theorem~\ref{thm:Tiba13MainThm}.
\end{proof}

\subsubsection[The case X=P\^{}3 and D=3H]{The case $\boldsymbol{X=\mathbb{P}^3}$ and $\boldsymbol{D=3H}$}
\label{subsubsec:P3-3H}

Let $X=\mathbb{P}^3$ and $D=3H$. A~general cubic has $20$ coefficients.
We record two complementary regimes: a~sparse support set $S$, for which
we verify the combinatorial hypotheses and obtain a~generic conclusion,
and the full support
$S=P_D\cap M$, where only the directions with $\ell_H=4$ can contribute
hypersurface components.

We first consider a~sparse support set. Let
\begin{align}
S={}&\{(0,0,0),(3,0,0),(0,3,0),(0,0,3),\nonumber\\
 &(1,0,0),(0,2,0),(2,0,1),(0,1,2)\}\subset P_D\cap M.\label{eq:S-P3}
\end{align}
Then $D_0$ is defined by
\begin{align}
	p(x,y,z,w)={}&a_{300}x^3+a_{030}y^3+a_{003}z^3+a_{000}w^3\nonumber\\
	&+a_{100}xw^2+a_{020}y^2w+a_{201}x^2z+a_{012}yz^2.\label{eq:p-P3-sparse}
\end{align}

\begin{Proposition}[a facewise criterion for the sparse cubic in $\mathbb{P}^3$]\label{prop:P3-sparse}
	Let $S$ be as in~\eqref{eq:S-P3} and let $p$ be as in~\eqref{eq:p-P3-sparse}.
	If, for every positive-dimensional face $\tau$ of $P_D$, the restricted
	coefficient vector lies outside $B_{\tau\cap S}$, then $
	T_N\setminus V(p)\hookrightarrow\mathbb{P}^3$
	is a~Kobayashi hyperbolic embedding.
In particular, this conclusion holds for a~nonempty Zariski open subset
of the coefficient torus where all eight coefficients are nonzero.
\end{Proposition}

\begin{proof}
	The support contains all four vertices of $3\Delta_3$, so every face
	is spanned by its restricted support and $\mathfrak d_S\ne\varnothing$.
	Normalize the nonzero cross products of pairs of difference vectors to
	primitive normals up to sign. Counting their projection layers gives
	\[
	\begin{array}{@{}ccc@{}}
\toprule
\text{support}&\text{layer count: number of directions}&
	\text{required minimum}\\ \midrule
	S&4:7,\ 5:24,\ 6:62&4\\
	S\cap\tau\ (\tau\text{ a facet})&3:2,\ 4:4\text{ per facet}&3\\
	S\cap\tau\ (\tau\text{ an edge})&\text{two or three points}&2\\
\bottomrule
\end{array}
	\]
	On each facet the two-dimensional enumeration uses primitive difference
	directions in its affine lattice. The four facets have the same layer
	count distribution. The $xy$- and $zw$-edges contain two support points
	each, and the other four edges contain three support points each.
	The supplementary script described in
	Appendix~\ref{app:verification} verifies these counts, including all
	$93$ directions for the eight-point support. Thus the facewise combinatorial
	hypotheses of Theorem~\ref{thm:Tiba13MainThm} hold.
	Lemma~\ref{lem:facewise-hyperbolic-embedding} proves the first assertion.
	The global exceptional locus $\mathcal Y_S\subset\mathbb P^7$ is proper
	and closed by its construction. Its complement meets the dense
	coefficient torus, giving the asserted nonempty Zariski open subset.
\end{proof}

Nonvanishing of the eight coefficients alone does not imply the facewise
bad-locus condition, because layer polynomials can cancel on the torus.
Indeed, set all eight coefficients in~\eqref{eq:p-P3-sparse} equal to $1$.
Then the translated one-dimensional subtorus
$
u\longmapsto[-u:u:-1:1]$, $ u\in\mathbb C^*$,
is contained in $V(p)$: direct substitution gives
$-u^3+u^3-1+1-u+u^2-u^2+u=0$. Thus the coefficient vector belongs to
the bad locus even though every coefficient is nonzero.
This example disproves the implication from coefficient nonvanishing to
bad-locus avoidance. It does not by itself prove that the complement is
nonhyperbolic. Accordingly, the proposition uses explicit facewise
bad-locus avoidance, and its generic conclusion does not assert the result
for every nonzero coefficient vector.

For the full support, the divisorial part is reduced to a~finite elimination
problem. When $S=P_D\cap M$, the arrangement $\mathcal H_S$ is large, but
the incidence estimate in Lemma~\ref{lem:YS-closed} shows that a~hyperplane
with $\ell_H\ge5$ contributes a~subset of codimension at least $2$; the
enumeration below shows that every nondivisorial direction actually has
$\ell_H\ge7$.
Consequently, only the finitely many hyperplanes with $\ell_H=4$ can
contribute hypersurface components.

\begin{Proposition}[divisorial degree for cubic surfaces]\label{prop:P3-divisorial-degree}
The codimension-one part of the global exceptional locus
$\mathcal Y\subset |3H|$ is reduced of degree $240$.
\end{Proposition}

\begin{proof}Let
	\[
	S=3\Delta_3\cap\mathbb Z^3
	=\bigl\{(i_1,i_2,i_3)\in\mathbb Z_{\ge0}^3\mid i_1+i_2+i_3\le3\bigr\}.
	\]
	Represent a~hyperplane $H\subset\mathbb R^3$ by a~primitive normal
	$n\in\mathbb Z^3$, up to sign. Then
	$\ell_H=|\{\langle n,I\rangle\mid I\in S\}|$. Direct enumeration gives
	exactly seven directions with $\ell_H=4$
$
	e_1$, $e_2$, $e_3$, $e_1+e_2$, $e_1+e_3$, $e_2+e_3$,
	$e_1+e_2+e_3$.
	For reproducibility, form all linearly independent pairs~${u,v\in(S-S)\setminus\{0\}}$, take $n=u\times v$, divide its coordinates
	by their greatest common divisor, and choose the sign so that the first
	nonzero coordinate is positive. Removing repetitions gives
	\[
	\begin{array}{@{}crrrrrrrrrr@{}}
\toprule
\ell_H&4&7&9&10&12&13&14&15&17&18\\ \midrule
	\#\{H\}&7&18&36&12&24&12&60&72&72&192.\\
\bottomrule
\end{array}
	\]
	Thus the literal definition of $\mathcal H_S$ gives $505$ hyperplane
	directions in total, exactly the seven displayed directions have four
	layers, and every other direction has at least seven layers.

	For one of these directions, write $P_0$, $P_1$, $P_2$, $P_3$ for the Newton
	polygons of the four layer polynomials in the rank-two lattice $H\cap M$.
	For the triples without a~point layer considered below, each support has
	dimension at least one and every pair has a~two-dimensional Minkowski sum.
	Thus the full triple is the essential subfamily. If its support
	differences generate $H\cap M$, the sparse resultant has multiplicity one
	and equals the reduced sparse eliminant
	\cite[Remark~3.5 and Proposition~3.6]{DAndreaJeronimoSombra2023}.
	In this case, the closure of $B_{S,H,j}$ is its irreducible zero locus,
	and the degree in each retained layer block is the mixed volume of the
	other two Newton polygons~\cite[equation~(3.2)]{DAndreaJeronimoSombra2023}.
 Here $\operatorname{MV}$ denotes the mixed volume with respect to
the lattice $H\cap M$, normalized by~${
\operatorname{MV}(P,P)=2\operatorname{Area}(P)}$,
where a~fundamental lattice parallelogram has area~$1$;
see~\cite[Section~2.2]{DAndreaJeronimoSombra2023}.
	We handle triples containing a~point layer separately, so that no
	resultant multiplicity is included in the reduced degree.

	For the four face-parallel directions
	$e_1$, $e_2$, $e_3$, $e_1+e_2+e_3$, the polygons are, up to reversing their order,
	$
	P_m=m\Delta_2$, $ 0\le m\le3$.
	If the point layer $P_0$ is omitted, the other three layers have degrees
	$1$, $2$, $3$. Their supports satisfy the preceding dimension conditions and
	generate the full lattice, since the degree-one support contains a~unimodular triangle. Their resultant therefore has total degree
	\[
	\operatorname{MV}(2\Delta_2,3\Delta_2)
	+\operatorname{MV}(\Delta_2,3\Delta_2)
	+\operatorname{MV}(\Delta_2,2\Delta_2)
	=6+3+2=11.
	\]
	If a~nonpoint layer is omitted, the point-layer polynomial must vanish.
	The other two retained polygons have positive mixed volume, so two general
	polynomials with those supports have a~common torus zero by
	Bernshtein's theorem~\cite[Theorem~A]{Bernshtein1975}. Consequently,
	after passing to the reduced closure, the three such cases give the same
	coordinate hyperplane, of degree $1$. Thus each face-parallel direction
	contributes $11+1=12$ rather than counting that coordinate hyperplane three
	times.

	For the other three directions, the polygons can be identified with
	$
	P_m=[0,3-m]\times[0,m]$, $ 0\le m\le3$.
	For rectangles $[0,a]\times[0,b]$ and $[0,c]\times[0,d]$, the normalized
	mixed volume is $ad+bc$. Hence
	\[
	\begin{array}{@{}ccccccc@{}}
\toprule
(a,b)&(0,1)&(0,2)&(0,3)&(1,2)&(1,3)&(2,3)\\ \midrule
	\operatorname{MV}(P_a,P_b)&3&6&9&5&6&3\\
\bottomrule
\end{array}
	\]
	and the four resultants obtained by omitting, respectively,
	$P_0$, $P_1$, $P_2$, $P_3$ have total degrees
	14, 18, 18, 14.
	For every omission, each retained polygon has dimension at least one and
	every retained pair has a~two-dimensional Minkowski sum. Moreover, at
	least one of $P_1$, $P_2$ is retained; its support contains differences equal
	to the two coordinate unit vectors. Thus each triple generates the full
	lattice and its resultant has multiplicity one. Each of these three
	directions consequently contributes
	$14+18+18+14=64$.

	It remains to distinguish the reduced components. Each noncoordinate
	eliminant depends on every individual retained coefficient. Indeed, fix
	the other two retained polynomials generically. Their positive mixed
	volume gives a~nonempty finite set of common torus zeros by
	Bernshtein's theorem~\cite[Theorem~A]{Bernshtein1975}. The Poisson
	formula~\cite[Theorem~3.3]{DAndreaJeronimoSombra2023} expresses the
	eliminant, up to a~nonzero factor independent of the remaining polynomial,
	as the product of its evaluations at those zeros, with multiplicities.
	Every monomial evaluates nontrivially there, so each of its coefficients
	occurs with a~nonzero leading coefficient in this product. Hence the
	coefficients absent from an eliminant are exactly those in its omitted layer.

	For the four noncoordinate face-parallel components, the omitted layers
	are the four different singleton vertex sets. For the twelve remaining
	components, they are six different four-point edge sets and six different
	six-point sets. Each six-point set consists of the cubic monomials of
	degree one in a~specified pair of homogeneous variables and degree two
	in the complementary pair. The variables that occur squared determine
	that complementary pair, so these six sets are distinct. Thus the
	sixteen noncoordinate components are pairwise distinct and are not
	coordinate hyperplanes. The four coordinate components are the four
	different vertex-coefficient hyperplanes; the repetitions within each
	face-parallel direction have already been counted only once.\looseness=-1

	The reduced divisorial part therefore has degree $
	4\cdot12+3\cdot64=240$.
	All remaining full-support hyperplanes have $\ell_H\ge7$. Since $r=3$,
	the incidence estimate in Lemma~\ref{lem:YS-closed} gives codimension at
	least $\ell_H-r\ge4$, so they cannot contribute divisors. On a~two-dimensional facet, the support is~${3\Delta_2\cap\mathbb Z^2}$ and its
	minimum layer count is $4$; the same estimate gives codimension at least~${4-2=2}$. On a~one-dimensional edge there are four lattice points, so its
	bad locus has codimension $4-1=3$. Finally, a~component of the coefficient
	boundary $\mathfrak d_S^*\setminus\mathfrak d_S$ requires all ten
	coefficients on one facet to vanish and therefore has codimension $10$.
	Thus neither a~proper-face bad locus nor the coefficient-chart boundary
	contributes a~divisor. This proves the assertion.
\end{proof}

The degree $240$ concerns the reduced divisorial part of the constructed
exceptional locus $\mathcal Y$. It is not the degree of an independently
identified locus of all nonhyperbolic complements. The computational
checks and representative lattice coordinates are given in
Appendix~\ref{app:verification} and the supplementary script.

% ------------------------------------------------------
\subsection{Hirzebruch surfaces}
\label{ex_hirzebruch}
\label{subsec:examples-hirzebruch}

Let $\mathbb{F}_\ell$ $(\ell\ge 1)$ be the Hirzebruch surface. We use the toric (GIT) presentation
\[
\mathbb{F}_\ell=\frac{\mathbb C^4\setminus\big(\mathbb C^2\times\{(0,0)\} \cup \{(0,0)\}\times\mathbb C^2\bigr)}
{\bigl\{\bigl(\lambda^\ell\mu,\mu,\lambda,\lambda\bigr)\mid \lambda,\mu\in\mathbb C^*\bigr\}}.
\]
The coordinate order is $(x,w,y,z)$.
There is a~toric contraction
$\pi\colon \mathbb F_\ell\to\mathbb P(\ell,1,1)$ of the negative section to
$[1:0:0]$. This is the unique singular point when $\ell\ge2$; for $\ell=1$
the target is~$\mathbb P^2$ and the image point is smooth.
On the chart $w\ne0$, the contraction is represented in weighted
homogeneous coordinates by
$
\pi(x,w,y,z)=[x/w:y:z]$.
This algebraic chart formula involves no choice of an $\ell$-th root of $w$.
The Picard group is generated by torus-invariant divisors; writing~${D_w=\{w=0\}}$ and $D_y=\{y=0\}$,
$
\Pic(\mathbb{F}_\ell)\cong \mathbb Z[D_w]\oplus \mathbb Z[D_y]$.
With the Cox grading
\[
\deg(x)=(1,\ell),\qquad \deg(w)=(1,0),\qquad
\deg(y)=\deg(z)=(0,1),
\]
the pair $(a,b)$ denotes the class $aD_w+bD_y$. Here
$D_w^2=-\ell$, $D_w\cdot D_y=1$, and $D_y^2=0$; consequently~${aD_w+bD_y}$ is ample if and only if $a>0$ and $b>\ell a$
(equivalently, for integral classes, $a\ge1$ and $b\ge\ell a+1$).
We concentrate on the ample class
$
D:=2D_w+(2\ell+2)D_y$,
for which $b-\ell a=2$.

We use the following six-term support set. Let
\[
S=\{(0,0),(0,1),(0,2),(2\ell+2,0),(\ell+2,1),(2,2)\},\]
where the coordinates of a~lattice point are the exponents $(i,j)$ of
the torus characters
\begin{equation}\label{eq:Hirzebruch-characters}
u:=y/z,\qquad v:=x/\bigl(wz^\ell\bigr).
\end{equation}
We consider the Cox-homogeneous polynomial
\begin{align*}
	p(x,w,y,z)
	=  a_{0,y}y^{2\ell+2}w^2\!+a_{0,z}z^{2\ell+2}w^2\!+a_{1,y}xy^{\ell+2}w+a_{1,z}xz^{\ell+2}w
	 +a_{2,y}x^2y^2+a_{2,z}x^2z^2.
\end{align*}
Every displayed monomial has Cox degree $(2,2\ell+2)$, so the polynomial
indeed defines a~divisor of class $D$.
Dividing by $w^2z^{2\ell+2}$ on the dense torus gives
\[
\frac{p}{w^2z^{2\ell+2}}
=a_{0,y}u^{2\ell+2}+a_{0,z}
+a_{1,y}u^{\ell+2}v+a_{1,z}v
+a_{2,y}u^2v^2+a_{2,z}v^2.
\]
More generally, a~Cox monomial of degree $(2,2\ell+2)$ has the form
$x^jw^{2-j}y^iz^{2\ell+2-\ell j-i}$ and corresponds to $u^iv^j$.
Consequently the polytope in these coordinates is
\begin{equation}\label{eq:Hirzebruch-polytope}
P_D=\{(i,j)\in\mathbb R^2\mid 0\le j\le2,\,
0\le i\le2\ell+2-\ell j\}.
\end{equation}
Its four vertices are $(0,0)$, $(2\ell+2,0)$, $(2,2)$, $(0,2)$, all of which
belong to $S$. For the direction determined by the three powers of $x$, the layer
polynomials are
\smash{$
a_{0,y}y^{2\ell+2}+a_{0,z}z^{2\ell+2}$},
\smash{$a_{1,y}y^{\ell+2}+a_{1,z}z^{\ell+2}$}, \smash{$
a_{2,y}y^2+a_{2,z}z^2$}
before specialization in the transverse variable. Thus one explicit part of
the facewise bad-locus avoidance condition is that, for nonzero coefficients,
no point $(y,z)\in\bigl(\mathbb C^*\bigr)^2$ makes two of these expressions vanish
simultaneously.

 A~toric boundary limit displays a~birational phenomenon. Set $w=1$ to
view the same equation on $\mathbb P(\ell,1,1)$.
The morphism $\pi$ is an isomorphism on the dense tori, identifying the
two dense torus complements. Their intrinsic hyperbolicity is therefore
the same. What depends on the compactification is the hyperbolic embedding
problem and the boundary criterion used to study it. The toric resolution
retains the exceptional negative section $w=0$, on which a~limiting entire
curve can appear. For $\ell=1$ this resolution is the blow-up of a~smooth
point of $\mathbb P^2$; for $\ell\ge2$ it resolves the singular point of
$\mathbb P(\ell,1,1)$.

Concretely, in the limiting configuration where the $x^2z^2$-term is
removed, the stability condition on the Cox quotient implies $x\ne0$ along
$w=0$, and the equation reduces to $a_{2,y}x^2y^2=0$. Outside
$\{y=0\}$ the divisor does not meet this boundary locus,
allowing the holomorphic map~${
\phi\colon \mathbb C\to\mathbb F_\ell}$, $ s\mapsto (x,w,y,z)=\bigl(1,0,1,e^s\bigr)$,
which lands in the boundary stratum away from the limiting divisor.
The curve $\phi$ violates the facewise ``no entire curves'' condition.
The proof of Proposition~\ref{prop:Hirzebruch-failure} establishes
failure of hyperbolic embedding of the torus complement in
$\mathbb F_\ell$.

\begin{Proposition}[a family losing hyperbolic embedding at a~toric boundary limit]
	\label{prop:Hirzebruch-failure}
	Fix nonzero complex numbers
	$a_{0,y}$, $a_{0,z}$, $a_{1,y}$, $a_{1,z}$, $a_{2,y}$ and consider the family
\begin{gather}
p_t=a_{0,y}y^{2\ell+2}w^2+a_{0,z}z^{2\ell+2}w^2
+a_{1,y}xy^{\ell+2}w+a_{1,z}xz^{\ell+2}w+a_{2,y}x^2y^2+t x^2z^2,\nonumber\\
D_t:=\divv(p_t).\label{eq:Hirzebruch-family}
\end{gather}
	Assume that there exists $t_0\in\mathbb C^*$ such that
	$D_{t_0}\in\mathfrak d_S^*\setminus\mathcal Y_S$, where
	$\mathcal Y_S$ is the facewise exceptional locus for the displayed
	six-point support. Then there is a~finite set
	$E\subset\mathbb C^*$ such that
	\[
	T_N\setminus\Supp(D_t\cap T_N)\hookrightarrow\mathbb F_\ell
	\]
	is a~Kobayashi hyperbolic embedding for every
$t\in\mathbb C^*\setminus E$. At $t=0$, the inclusion
\[
U_0:=T_N\setminus\operatorname{Supp}(D_0\cap T_N)
\hookrightarrow\mathbb F_\ell
\]
is not a~Kobayashi hyperbolic embedding for every choice of the five
nonzero coefficients; the failure at $t=0$ does not require the
assumption concerning $t_0$.
\end{Proposition}

\begin{proof}
We first verify generic hyperbolic embedding for $t\ne0$ in
\eqref{eq:Hirzebruch-family}.
	For $\ell\ge1$, the support~$S$ spans the polygon
	$P_D$ in~\eqref{eq:Hirzebruch-polytope}. For a~nonzero linear functional $\nu(i,j)=\alpha i+\beta j$,
	if $\beta\ne0$, the points $(0,0)$, $(0,1)$, $(0,2)$ already give three
	distinct values. If $\beta=0$, the four distinct first coordinates
	0, 2, $\ell+2$, $2\ell+2$ give four values. The quotient $\nu(i,j)=j$ has
	exactly three layers and its kernel belongs to $\mathcal H_S$. Hence
	$\min_{H\in\mathcal H_S}|\phi_H(S)|=3$ for every $\ell\ge1$.
The support $S$ contains all four vertices of $P_D$ in
\eqref{eq:Hirzebruch-polytope}.
	Thus each edge support contains both endpoints and has at least two
	layers; in particular, $\mathfrak d_S\ne\varnothing$.
	Thus the ordinary facewise hypotheses hold and~${\mathcal Y_S\subset\mathfrak d_S^*}$ is proper.

	The affine family extends to the morphism
	\[
	\overline\gamma\colon\ \mathbb P^1\longrightarrow \mathfrak d_S^*,
	\qquad
	[s:t]\longmapsto
	[sa_{0,y}:sa_{0,z}:sa_{1,y}:sa_{1,z}:sa_{2,y}:t],
	\]
	where the six entries are the coefficients in the displayed order. On the chart
	$s=1$, this is the divisor $D_t$. Since
	$D_{t_0}\notin\mathcal Y_S$, the image of $\overline\gamma$ is not contained
	in the closed set $\mathcal Y_S$. Therefore,~$\overline\gamma^{-1}(\mathcal Y_S)$ is a~proper closed subset of
	$\mathbb P^1$ and hence is finite. Set $
	E:=\{t\in\mathbb C^*\mid D_t\in\mathcal Y_S\}$.
	For $t\in\mathbb C^*\setminus E$, the construction of $\mathcal Y_S$ and
	Lemma~\ref{lem:facewise-hyperbolic-embedding} give the asserted embedding.

We now prove failure of hyperbolic embedding at $t=0$ directly from
the neighborhood criterion of
Section~\ref{subsec:hyperbolic-embeddings}. Fix $R>1$ and put
\smash{$
A_R:=|a_{1,y}|+|a_{1,z}|{\rm e}^{(\ell+2)R}$}, \smash{$
B_R:=|a_{0,y}|+|a_{0,z}|{\rm e}^{(2\ell+2)R}$}.
Choose $\varepsilon_R\ne0$ so that
\[
|\varepsilon_R|<\min\left\{1,\frac1R,
\frac{|a_{2,y}|}{2(1+A_R+B_R)}\right\}.
\]
For $|s|\le R$, direct substitution gives
\begin{align*}
p_0\bigl(1,\varepsilon_R,1,e^s\bigr)
={}&a_{2,y}
+\varepsilon_R\bigl(a_{1,y}+a_{1,z}{\rm e}^{(\ell+2)s}\bigr)+\varepsilon_R^2\bigl(a_{0,y}+a_{0,z}{\rm e}^{(2\ell+2)s}\bigr).
\end{align*}
The choice of $\varepsilon_R$ gives
\[
|p_0\bigl(1,\varepsilon_R,1,e^s\bigr)-a_{2,y}|
\le |\varepsilon_R|A_R+|\varepsilon_R|^2B_R
<|a_{2,y}|/2,
 \qquad |s|\le R.
\]
Consequently
$
f_R\colon \{|s|<R\}\longrightarrow U_0$, $
f_R(s)=\bigl[\bigl(1,\varepsilon_R,1,e^s\bigr)\bigr]$
is well defined. Here brackets denote the Cox quotient, and all four
coordinates of $f_R(s)$ are nonzero. As $R\to\infty$,
$
f_R(0)\longrightarrow[(1,0,1,1)]$, $
f_R(1)\longrightarrow[(1,0,1,e)]$.
The boundary points $[(1,0,1,1)]$ and $[(1,0,1,e)]$ are distinct
because their~$z/y$ coordinates are $1$ and $e$.
The convergences place both points in $\overline{U_0}$.
Distance decrease for the disk map~${\zeta\mapsto f_R(R\zeta)}$ yields
$
d_{U_0}(f_R(0),f_R(1))
\le\rho_{\mathbb D}(0,1/R)\longrightarrow0$.
For any neighborhoods $V_0$, $V_1$ of $[(1,0,1,1)]$ and~${[(1,0,1,e)]}$,
respectively, $f_R(0)\in V_0$ and $f_R(1)\in V_1$ for sufficiently
large~$R$. Hence
$
d_{U_0}(V_0\cap U_0,V_1\cap U_0)=0$,
contrary to the positive-distance condition in
Section~\ref{subsec:hyperbolic-embeddings}.
\end{proof}

The contraction sends the limiting curve $\phi$ to $[1:0:0]$.
Write $p_t$ for the displayed Cox polynomial and
$q_t(X,Y,Z):=p_t(X,1,Y,Z)$. Every monomial of $q_t$ contains $Y$ or
$Z$, so $[1:0:0]$ belongs to the closure of the torus divisor for every
$t$. In particular, $\pi$ does not map the whole complement
$\mathbb F_\ell\setminus\operatorname{Supp}D_0$ into the complement of
the pushforward divisor. It does identify the dense torus complements.
Proposition~\ref{prop:Hirzebruch-contraction} compares the embeddings
of $U_0$ in the two compactifications.

\begin{Proposition}[comparison after contraction]
\label{prop:Hirzebruch-contraction}
For each $\ell\ge1$, there is a~nonempty Zariski open subset
\smash{$
\Omega\subset\bigl(\mathbb C^*\bigr)^5
$}
of coefficient tuples $(a_{0,y},a_{0,z},a_{1,y},a_{1,z},a_{2,y})$
such that the hypothesis concerning $t_0$ in
Proposition~{\rm\ref{prop:Hirzebruch-failure}} holds and, at $t=0$,
$
U_0\hookrightarrow\mathbb P(\ell,1,1)
$
is a~Kobayashi hyperbolic embedding.
For the same coefficients, $U_0\hookrightarrow\mathbb F_\ell$ is not a~hyperbolic embedding. Thus the dense torus complement is the same
complex manifold in the two compactifications, but its hyperbolic
embedding behavior differs.
\end{Proposition}

\begin{proof}
We first exclude entire curves in the dense-torus divisor and its
complement for general coefficients. At $t=0$ its support in the
$(u,v)$ characters from~\eqref{eq:Hirzebruch-characters} is
\[
S_0=S\setminus\{(0,2)\}
=\{(0,0),(0,1),(2\ell+2,0),(\ell+2,1),(2,2)\}.
\]
Every nonzero linear functional $\nu(i,j)=\alpha i+\beta j$ takes at
least three distinct values on $S_0$. If~${\alpha=0}$, use the three
distinct second coordinates $0$, $1$, $2$; if $\beta=0$, use the distinct
first coordinates 0, $\ell+2$, $2\ell+2$.
If $\alpha\beta\ne0$, the values at (0,0), (0,1), $(2\ell+2,0)$ are
distinct unless $\beta=(2\ell+2)\alpha$. In that case the value at
$(\ell+2,1)$ is $(3\ell+4)\alpha$, distinct from the first two values.
The support has affine dimension two, so
Lemma~\ref{lem:YS-closed} gives a~proper closed local exceptional locus
$Y_{S_0}\subset\mathbb P^4$. Let~$\Omega_2$ be the nonempty Zariski
open subset of the coefficient torus on which $[a]\notin Y_{S_0}$.
Lemma~\ref{lem:no-entire-curves-in-torus} then excludes entire curves
both in the torus divisor and in $U_0$.

We next verify the boundary conditions in the contracted surface.
Its torus-divisor closure is defined by the weighted homogeneous equation
\begin{align*}
q_0(X,Y,Z)={}&a_{0,y}Y^{2\ell+2}+a_{0,z}Z^{2\ell+2}
+a_{1,y}XY^{\ell+2}+a_{1,z}XZ^{\ell+2}+a_{2,y}X^2Y^2.
\end{align*}
No coordinate divisor is a~component, since none of $X$, $Y$, $Z$ divides
$q_0$. On the one-dimensional orbits in $X=0$, $Y=0$, and $Z=0$,
respectively, the equations reduce, after an invertible monomial factor,
to
\begin{gather*}
a_{0,y}u^{2\ell+2}+a_{0,z}, \qquad u=Y/Z,\qquad
a_{0,z}+a_{1,z}\xi,\qquad\xi=X/Z^\ell,\\
a_{0,y}+a_{1,y}\eta+a_{2,y}\eta^2,
\qquad \eta=X/Y^\ell.
\end{gather*}
All these coefficients are nonzero. Each zero set in $\mathbb C^*$
is consequently finite and nonempty. It contains no nonconstant entire
curve, and its complement in $\mathbb C^*$ is a~sphere with at least
three punctures, hence is Brody hyperbolic. The remaining orbits are
points. Tiba's boundary criterion~\cite[Lemma~3]{TibaYusaku13}, applied
to the reduced torus divisor and its closure in the projective toric
surface $\mathbb P(\ell,1,1)$, gives the asserted embedding for
coefficients in $\Omega_2$. The hypotheses of
\cite[Lemma~3]{TibaYusaku13} concern the reduced torus divisor and its
closure, and impose no Cartier assumption on the weighted-projective
divisor class.

Finally, we arrange the $t_0$ hypothesis simultaneously. Choose finitely
many homogeneous equations $G_j$ defining the proper global locus
$\mathcal Y_S\subset\mathbb P^5$ for the six-point support.
Define
\[
\Omega_1:=\{a\in\bigl(\mathbb C^*\bigr)^5\mid
G_j(a,t)\not\equiv0\text{ for some }j\}.
\]
The set $\Omega_1$ is Zariski open because the condition is the
nonvanishing of at least one coefficient of a~polynomial $G_j(a,t)$
viewed as a~polynomial in $t$. It is nonempty because the
complement of the proper closed set $\mathcal Y_S$ meets the dense
six-coordinate coefficient torus. For $a\in\Omega_1$, only finitely
many $t$ belong to $\mathcal Y_S$, so some $t_0\in\mathbb C^*$ avoids
it. The coefficient torus is irreducible, and hence~${\Omega:=\Omega_1\cap\Omega_2}$ is nonempty and Zariski open.
Proposition~\ref{prop:Hirzebruch-failure} supplies the failure of
embedding in~$\mathbb F_\ell$ at $t=0$ for these same tuples.
\end{proof}

For the Cartier comparison, we also record the auxiliary ample class
$
D-D_y=2D_w+(2\ell+1)D_y$.
It is included only to compare the Cartier behavior after decreasing the
fiber coefficient by one. Every member $D_t$ of the displayed pencil,
including $D_0$, still has class $D$; the class $D-D_y$ is not produced by
setting $t=0$. In the following table the symbols $\mathcal O(d)$ denote
the rank-one reflexive sheaves associated with the pushforward Weil
classes. Since $D_w$ is contracted and $D_y$ pushes forward to a~weight-one coordinate divisor, the two pushforward degrees are
$2\ell+2$ and $2\ell+1$, respectively. On the well-formed weighted
projective surface $\mathbb P(\ell,1,1)$, the rank-one reflexive sheaf
$\mathcal O(d)$ is invertible precisely when $\ell\mid d$.
Indeed, at $[1:0:0]$ the local cyclic quotient has character group
$\mathbb Z/\ell\mathbb Z$, and the class of $\mathcal O(d)$ is $d$
modulo $\ell$; away from this point the surface is smooth.
The pushforwards of the two classes therefore satisfy
\[
\begin{array}{@{}ccc@{}}
\toprule
\text{class on }\mathbb F_\ell&\text{pushforward class}&
\begin{gathered}\text{values of }\ell\ge1\text{ for which}\\\text{the class is Cartier}\end{gathered}\\ \midrule
2D_w+(2\ell+2)D_y&\mathcal O(2\ell+2)&\ell=1,2\\
2D_w+(2\ell+1)D_y&\mathcal O(2\ell+1)&\ell=1\\
\bottomrule
\end{array}
\]
For $\ell\ge3$, even the first pushforward class is already non-Cartier.
Moreover, the ample divisor on $\mathbb F_\ell$ is not the pullback of
the corresponding class on the contracted surface: for $E=D_w$, its
intersection with $E$ is $2$ (respectively $1$ for the second class),
whereas a~Cartier pullback has intersection zero with $E$.
Consequently, the integral-polytope criterion for ample Cartier classes
cannot be transferred directly to these non-Cartier weighted-projective
classes. Proposition~\ref{prop:Hirzebruch-contraction} instead uses the
local Laurent-polynomial criterion and Tiba's boundary criterion for the
closure of the torus divisor, with the boundary orbits checked explicitly.
	
\appendix
\section{Reproducibility of the finite computations}
\label{app:verification}

The supplementary file \texttt{verify\_toric\_computations.py} uses only
the Python standard library (Python 3.9 or later) and exact integer
arithmetic. Running
\begin{verbatim}
python3 verify_toric_computations.py --output verification_results.json
\end{verbatim}
checks the assertions below and writes the direction list, layer supports,
lattice coordinates, and mixed-volume data to the indicated JSON file.
The supplied \texttt{verification\_results.json} is the output of this
command. The program checks the finite combinatorics and the degree
bookkeeping; it does not expand symbolic resultants or replace the analytic
proofs and the sparse-resultant arguments in the text.

For a~three-dimensional support, the program forms all independent pairs
in $S-S$, normalizes their cross products to primitive normals up to sign,
and removes duplicates. It then groups $S$ by the integer values
$\langle n,I\rangle$. In dimension two it uses the primitive normals to
the nonzero difference vectors. These procedures enumerate exactly the
hyperplanes in the definition of $\mathcal H_S$. They verify the $505$
directions for $3\Delta_3\cap\mathbb Z^3$, the $93$ directions for the
eight-point support, the four sparse facets and six edges, and each stage
of the five-coefficient deletion sequence in
Proposition~\ref{prop:P2-3H-five-term}.

Here are lattice coordinates for the two types of four-layer direction
in Proposition~\ref{prop:P3-divisorial-degree}. For~${n=e_3}$, the layer
$i_3=3-m$ has coordinates $(i_1,i_2)$ and polygon $m\Delta_2$.
For $n=e_1+e_2$, the layer $i_1+i_2=m$ is based at $(0,m,0)$.
In the basis $(0,0,1)$, $(1,-1,0)$ of $H\cap M$, its coordinates are~$(i_3,i_1)$ and its polygon is $[0,3-m]\times[0,m]$.
Permutations of the four homogeneous coordinates give the other five
directions. Alternatively, for any of the seven primitive normals $n$
(all of whose entries are $0$ or $1$), choose $p$ with $n_p=1$.
The vectors $e_j-n_j e_p$, $j\ne p$, form a~basis of~${H\cap M}$.
The program uses this basis separately for every direction.

For each triple of retained layer supports with no point layer, the
program checks the essential-subfamily dimension conditions and that the
greatest common divisor of the two-by-two determinants of their within-layer
differences is $1$. The latter verifies the full-lattice condition in the
resultant multiplicity formula. Polygon areas are computed by the exact
shoelace formula after taking convex hulls, and
$
\operatorname{MV}(P,Q)
=\operatorname{Area}(P+Q)-\operatorname{Area}(P)-\operatorname{Area}(Q)$.
The resulting reduced components have the following degrees:
\[
\begin{array}{@{}crrr@{}}
\toprule
\text{component type}&\text{number}&\text{degree each}&\text{contribution}\\ \midrule
\text{coordinate hyperplane}&4&1&4\\
\text{face-parallel noncoordinate eliminant}&4&11&44\\
\text{other noncoordinate eliminant}&6&14&84\\
\text{other noncoordinate eliminant}&6&18&108\\ \midrule
\text{total}&20&&240\\
\bottomrule
\end{array}
\]
Repeated coordinate components are counted once. The omitted coefficient
sets of the sixteen noncoordinate components are checked to be distinct;
their use to distinguish the components is justified in the proof of
Proposition~\ref{prop:P3-divisorial-degree}. The program also verifies the
layer bounds giving codimension at least $4$ for other full-support
directions, $2$ on facets, and $3$ on edges. The coefficient-chart boundary
has codimension $10$. Thus these additional sets do not change the
divisorial degree.

\subsection*{Acknowledgements}
The authors thank the anonymous referees and the SIGMA editor for their
careful reading and helpful comments, which improved the presentation of
this paper.

%\subsection*{Conflicts of interest}
%On behalf of all authors, the corresponding author states that there is no conflict of interest.

\pdfbookmark[1]{References}{references}
\LastPageEnding
	
\end{document}